\documentclass[12pt,twoside]{article}
\usepackage{amsmath}
\usepackage{amsfonts}
\usepackage{amsthm}
\usepackage{amssymb}
\usepackage{multicol}

\begin{document}
\renewcommand{\baselinestretch}{1.5}
\title{{\normalsize
{\bf Splitting criteria for a definite 4-manifold with infinite cyclic 
fundamental group\footnote{This paper is dedicated to the memory of 
Dr. Tim D. Cochran.}}}}
\author{{\footnotesize Akio KAWAUCHI}\\
\date{}
{\footnotesize{\it Osaka City University Advanced Mathematical Institute}}\\
{\footnotesize{\it Sugimoto, Sumiyoshi-ku, Osaka 558-8585, Japan}}\\
{\footnotesize{\it kawauchi@sci.osaka-cu.ac.jp}}}
\maketitle
\vspace{0.25in}
\baselineskip=15pt
\thispagestyle{empty}
\newtheorem{Theorem}{Theorem}[section]
\newtheorem{Conjecture}[Theorem]{Conjecture}
\newtheorem{Lemma}[Theorem]{Lemma}
\newtheorem{Sublemma}[Theorem]{Sublemma}
\newtheorem{Proposition}[Theorem]{Proposition}
\newtheorem{Corollary}[Theorem]{Corollary}
\newtheorem{Claim}[Theorem]{Claim}
\newtheorem{Definition}[Theorem]{Definition}
\newtheorem{Example}[Theorem]{Example}
\begin{abstract}  
Two criteria for a closed connected definite 
4-manifold with infinite cyclic fundamental group to be TOP-split are  
given. One criterion extends a sufficient condition made in a previous paper. 
The result is equivalent to a purely algebraic result 
on the question asking when a positive definite Hermitian form
over the ring of integral one-variable Laurent polynomials is represented 
by an integer matrix. As an application, an infinite family of 
orthogonally indecomposable unimodular odd definite symmetric $Z$-forms is produced. 
\end{abstract}

\phantom{x}

\phantom{x}


\phantom{x}


\phantom{x}

\phantom{x}

\noindent{\bf 1. Introduction}

A closed connected oriented topological 4-manifold $M$ is called a 
$Z^{H_1}$-{\it manifold} if there is a fixed isomorphism from the first homology group 
$H_1(M)$ onto $Z$, and a 
$Z^{\pi_1}$-{\it manifold} if there is a fixed isomorphism from the fundamental group 
$\pi_1(M)$ onto $Z$. 
A $Z^{\pi_1}$-manifold $M$ is {\it TOP-split} if $M$ is homeomorphic to the 
connected sum $S^1\times S^3\#M_1$ for 
a simply connected closed 4-manifold $M_1$  obtained from $M$ by 
a $2$-handle surgery killing $\pi_1(M)\cong Z$, and {\it virtually} TOP-split if 
a finite connected  covering space 
of  $M$ is TOP-split.  A $Z^{H_1}$-manifold $M$  is {\it definite} 
if the rank of the $Z$-intersection form 
\[\mbox{Int}^M:H_2(M;Z)\times H_2(M;Z)\to Z\]
is equal to the absolute value of the signature, 
and {\it positive definite} if, furthermore, the signature is positive. 
A definite $Z^{\pi_1}$-manifold with negative signature is changed to be  positive definite 
by reversing the  orientation of $M$. 

Before explaining the main theorem of this paper, 
a history on the TOP-splitting of a $Z^{\pi_1}$-manifold is described here. 

In \cite{K0}, every topological $Z^{\pi_1}$-manifold was claimed to be TOP-split. However, 
this is not true, as noted in \cite{KK,K1}.  Concerning this error, 
I. Hambleton and P. Teichner in  \cite{HT} have constructed an example of a $\Lambda$-Hermitian matrix $L$ 
with determinant $+1$  
which cannot be $\Lambda$-conjugate to an integral  matrix, where $\Lambda=Z[Z]=Z[t,t^{-1}]$ denotes the integral Laurent polynomial ring. 
By  a construction  of  M.  H. Freedman and  F. Quinn in  \cite{FQ}, every$\Lambda$-Hermitian matrix $A$ is realized by 
a $\Lambda$-intersection matrix on  a  unique (up to  Kirby-Siebenmann 
obstructions $Z^{\pi_1}$-manifold $M_A$.  This 
$Z^{\pi_1}$-manifold $M_L$ is referred to as 
the {\it Hambleton-Teichner-Freedman-Quinn} $Z^{\pi_1}$-{\it manifold}.
Since every $\Lambda$-intersection matrix on  a  TOP-split $Z^{\pi_1}$-manifold is $\Lambda$-conjugate to an integral  matrix, 
the  $Z^{\pi_1}$-manifold $M_L$ is not TOP-split, and thus gives a counterexample 
to the TOP-splitting claim of a topological $Z^{\pi_1}$-manifold. 
Furthermore, the  $Z^{\pi_1}$-manifold $M_L$ was a positive 
definite, non-smoothable and virtually non-TOP-split $Z^{\pi_1}$-manifold,  which is shown 
 by S.  Friedl, I. Hambleton, P. Melvin,  and P. Teichner in \cite{FHMT}. 
In \cite{KA}, it was shown that every $Z^{\pi_1}$-manifold is TOP-split 
if and only if it is virtually TOP-split, which implies that 
every indefinite $Z^{\pi_1}$-manifold is TOP-split. Further, in \cite{KA2},  
a positive definite $Z^{\pi_1}$-manifold is TOP-split if  every finite covering space of it  admits an intersection matrix whose 
diagonal entries are smaller than or equal to $2$.  
As a consequence of these results in \cite{KA,KA2}, it was shown that every {\it smooth} 
$Z^{\pi_1}$-manifold is TOP-split. 

In this paper, the necessary and sufficient conditions on the TOP-splitting 
for a positive definite $Z^{\pi_1}$-manifold generalizing the sufficient condition of \cite{KA2} are given. 

For this purpose, the following two notions are useful: 

\medskip 

\noindent{$\bullet$} One notion is a notion of  a {\it winding degree} on a positive definite $Z^{\pi_1}$-manifold $M$  
which is a non-negative integer $\lambda$ measuring  a difference between $M$ and  
the TOP-split $Z^{\pi_1}$-manifold $S^1\times S^3\#M_1$. 
The winding degree $\lambda$ is defined in \S2 to take non-unique value  for a given $M$ for convenience,  but   
by definition the  minimum $\lambda_{\mbox{min}}$ of  all winding degrees $\lambda$ on $M$ is seen to be 
an invariant of $M$.  
It will be shown in Theorem~1.1 that $\lambda_{\mbox{min}}=0$ on $M$  if and only if  $M$ is TOP-split.  

\medskip 

\noindent{$\bullet$} The other notion is a notion of  a homology class called a {\it minimal element} in the second 
homology group $H_2(M;Z)$ of a positive definite $Z^{\pi_1}$-manifold $M$. This notion is   
a standard notion used for the proof of 
Eichler's unique orthogonal indecomposable splitting theorem for 
a positive definite symmetric bilinear form   
(see Eichler \cite{E},  Kneser  \cite{Kn}, Milnor-Husemoller \cite{MH}). 

\phantom{x}

\noindent{\bf Definition.} For a  positive definite $Z^{H_1}$-manifold $M$,  the definition of  a minimal element 
is given by the following two notions: 

\medskip

\noindent{$\bullet$} The {\it square length} of an element $x\in H_2(M;Z)$,  
denoted by $||x||^2$,  is 
the $Z$-self-intersection number $\mbox{Int}^M(x,x)$. 

\noindent{$\bullet$} An element $x\in H_2(M;Z)$ is {\it minimal} if $x\ne 0$ and  $x$ cannot be 
the sum $y+z$ of any elements $y,z\in H_2(M;Z)$ such that 
\[||x||^2>||y||^2\quad {\mbox{and}} \quad ||x||^2>||z||^2.\] 

\phantom{x}

As a basic observation, every minimal element of $H_2(M;Z)$ belongs to the unique indecomposable orthogonal sum component of 
$H_2(M;Z)$.

 \phantom{x}
 
\noindent{\bf Notation.} For a $Z^{H_1}$-manifold $M$,  the following notations are used.

\medskip

\noindent{$\bullet$} The $m$-fold cyclic  connected covering space of $M$ is denoted by $M^{(m)}$.  

\noindent{$\bullet$} The infinite cyclic connected covering space of $M$ with covering transformation 
group generated by $t$ is denoted by  $\widetilde M$. 

\phantom{x}

For every $Z^{\pi_1}$-manifold $M$, it is known that  the $\Lambda$-module $H_2(\widetilde M;Z)$    is a free 
$\Lambda$-module 
\[H_2(\widetilde M;Z)\cong \Lambda^n\]
of rank $n=\beta_2(M)$.
This fact was proved  in  \cite[Lemma~2.1]{Kold2} for a more general oriented compact 4-manifold with infinite cyclic 
fundamental group by using  three integral dualities on an infinite cyclic covering of  a topological manifold in \cite{Kold1}.  
Also, see  \cite{FQ} for another proof. 

In \cite{KA}, it is shown (as stated above) that 
$M$ is TOP-split if and only if $M^{(m)}$ is TOP-split 
for some $m$. 
Let 
\[\Lambda^{(m)}=\Lambda/(t^m-1)\Lambda\] 
be the quotient ring 
of the Laurent polynomial ring $\Lambda$ by the ideal $(t^m-1)\Lambda$.  
For a $Z^{\pi_1}$-manifold $M$, the $\Lambda$-module $H_2(M^{(m)};Z)$ is 
identical to the quotient $\Lambda$-module 
\[H_2(\widetilde M;Z)/(t^m-1)H_2(\widetilde M;Z), \]
which is a free $\Lambda^{(m)}$-module of rank $n$.  
For an element $\tilde x\in H_2(\widetilde M;Z)$,  let $\tilde x^{(m)}\in H_2(M^{(m)};Z)$ denote 
the projection image of $\tilde x$ under the covering projection homomorphism 
$H_2(\widetilde M;Z)\to H_2(M^{(m)};Z)$.  

\phantom{x}

\noindent{\bf Definition.} For a positive definite $Z^{\pi_1}$-manifold $M$,  
the infinite cyclic covering space $\widetilde M$ of $M$ and  $\Lambda=Z[t,t^{-1}]$,  the following definitions 
are set. 

\medskip

\noindent{$\bullet$} The $\Lambda$-{\it square length} of an element $\tilde x\in H_2(\widetilde M;Z)$  
denoted by $||\tilde x||^2_{\Lambda}$  is the $\Lambda$-self-intersection number 
$\mbox{Int}^{\widetilde M}_{\Lambda}(\tilde x,\tilde x)$, which  
is an integral Laurent polynomial $a(t)\in \Lambda$ in $t$ with $t$-symmetry 
$a(t)=a(t^{-1})$ (see \S~2). 
 
\noindent{$\bullet$} The {\it exponent} of  $\tilde x$, denoted by $e(\tilde x)$,  
is the highest degree of $a(t)=||\tilde x||^2_{\Lambda}$ which is a non-negative integer. 

\noindent{$\bullet$} An element $\tilde x\in H_2(\widetilde M;Z)$  is {\it minimal} if $\tilde x\ne 0$ and 
$\tilde x$ cannot be the sum $\tilde y+ \tilde z$ 
of any elements $\tilde y, \tilde z\in H_2(\widetilde M;Z)$ such that 
\[||\tilde x||^2>||\tilde y||^2 \quad {\mbox{and}} \quad ||\tilde x||^2>||\tilde z||^2.\] 

\phantom{x}

It is noted in \S~2 that 
the notion of a minimal element in $H_2(\widetilde M;Z)$ is a natural generalization 
of the notion of a minimal element in $H_2(M;Z)$, for a positive definite $Z^{\pi_1}$-manifold $M$.  
By the positivity of a square length shown in \S~2, it is seen that 
every non-zero element $\tilde x\in H_2(\widetilde M;Z)$ is the sum 
of finitely many minimal elements and every minimal element of $H_2(\widetilde M;Z)$ 
belongs to a unique indecomposable orthogonal sum component of $H_2(\widetilde M;Z)$.  
The multiplication $t^k \tilde x$ 
for an integer $k$ is called a $t$-{\it power shift} of $\tilde x$. 
The double covering projection $p:M^{(2m)}\to M^{(m)}$ is particularly 
used in the arguments of this paper. 
We shall show the following theorem:

\phantom{x}

\noindent{\bf Theorem~1.1.} 
The following conditions (0)-(5) 
on a positive definite $Z^{\pi_1}$-manifold $M$ are mutually equivalent:

\medskip

\noindent(0) The $Z^{\pi_1}$-manifold $M$ is TOP-split.

\noindent(1) There are elements
$\tilde x_i\, (i=1,2,\dots,n)$ in  $H_2(\widetilde M;Z)$ such that 
the elements 
$\tilde x_i^{(1)}\, (i=1,2,\dots,n)$ are $Z$-generators 
for  $H_2(M;Z)$ and the $\Lambda$-intersection numbers 
$\mbox{Int}^{\widetilde M}_{\Lambda}(\tilde x_i, \tilde x_j)$ 
are integers for all $i,j$. 

\noindent(2) Any minimal $\Lambda$-generators 
$\tilde x_i\, (i=1,2,\dots,n)$ of  $H_2(\widetilde M;Z)$ have the 
property that  after suitable $t$-power shifts of 
$\tilde x_i\, (i=1,2,\dots,n)$, 
the $\Lambda$-intersection numbers 
$\mbox{Int}^{\widetilde M}_{\Lambda}(\tilde x_i, \tilde x_j)$ 
are integers for all $i,j$.   

\noindent(3) Every minimal element 
$\tilde x$ of  $H_2(\widetilde M;Z)$ is sent to a minimal element 
$\tilde x^{(m)}\in H_2(M^{(m)};Z)$  
for every $m$ such that every element $x'\in H_2(M^{(2m)};Z)$ 
with $p_*(x')=\tilde x^{(m)}$ satisfies the inequality 
\[||x'||^2 \geqq ||\tilde x^{(m)}||^2.\]

\noindent(4) For  any previously given winding degree $\lambda$ on $M$,  
there is an $m \geqq \lambda$ for which there are  minimal 
$Z$-generators  $x_i\, (i=1,2,\dots,s)$ of  
$H_2(M^{(m)};Z)$ such that  for every $i$, and for every element 
$x'_i \in H_2(M^{(2m)};Z)$ with $p_*(x'_i)=x_i$,  the inequality 
\[||x'_i||^2 > ||x_i||^2-2\]
holds. 

\noindent(5) The minimal winding degree $\lambda_{\mbox{min}}$ on $M$ is zero. 

\phantom{x}

In Theorem~1.1,  (3), (4) and (5) are new results. Note that (0) is equivalent to (1) 
without assumption of the positive definiteness on $M$ and  is equivalent to the following condition:

\phantom{x}

(1)${}^*$ There is a $\Lambda$-basis $\tilde x_i\, (i=1,2,\dots,n)$ of  $H_2(\widetilde M;Z)$ such that 
the $\Lambda$-intersection numbers 
$\mbox{Int}^{\widetilde M}_{\Lambda}(\tilde x_i, \tilde x_j)$ 
are integers for all $i,j$. 

\phantom{x}

It is known in \cite{K1,K2} that a $Z^{\pi_1}$-manifold $M$ with property (1)${}^*$ is TOP-split (see the proof of 
(1) $\to$ (0) in the proof of Theorem~1.1). 
The $\Lambda$-intersection form on $M$ with with property (1)${}^*$ is said to be $Z$-{\it extended} in \cite{HT}, 
and also said to be  {\it exact}  because such a manifold is a special case of a $Z^{H_1}$-manifold admitting  
an exact sequence called an exact $Z^{H_1}$-manifold, which is discussed in \cite{K2}. 

It follows directly from (2) that for a positive definite $Z^{\pi_1}$-manifold $M$,  
the connected sum $Y\# M$ for any closed simply connected 
positive definite 4-manifold $Y$ is TOP-split if and only if $M$ is TOP-split. 

For a positive definite $Z^{\pi_1}$-manifold $M$, note that, by definition, if there are  
$Z$-generators $x_i\, (i=1,2,\dots, n)$ of $H_2(M^{(m)};Z)$ with 
$||x_i||^2\leqq 2$ for all $i$, then there are minimal $Z$-generators $y_j\, (j=1,2,\dots, s)$ of 
$H_2(M^{(m)};Z)$ with $||y_j||^2\leqq 2$ for all $j$. Then  for every element $y'_j\in H_2(M^{(2m)};Z)$ 
with $p_*(y'_j)=y_j$, the inequality 
\[||y'_j||^2 > 0 \geqq ||y_j||^2-2\]
holds. 
It will be explained in \S2 that a winding degree $\lambda$ on $M$ is taken 
smaller than or equal to a winding index $\delta$ on $M$ defined in \cite{KA2}.  
With these observations,  
the following corollary meaning the main theorem of \cite{KA2} is obtained as a consequence of Theorem~1.1 (4):

\phantom{x}

\noindent{\bf Corollary~1.2 (\cite[Theorem~1.1]{KA2}).} A positive 
definite $Z^{\pi_1}$-manifold $M$ is TOP-split if 
for any previously given 
winding index $\delta$ on $M$, there is  an $m\geqq \delta$  for which 
there is a  $Z$-basis $x_i$ 
$(i=1,2,\dots, n)$ of  $H_2(M^{(m)};Z)$  such that 
$||x_i||^2 \leqq 2$ for all $i$. 

\phantom{x}

By using Corollary~1.2, it was shown in  \cite{KA2} that 
every $Z^{\pi_1}$-manifold $M$ is TOP-split  
if for every $m$ the intersection form of $M^{(m)}$ 
is represented by a block sum of copies of $(1)$ and/or $E_8$. This means that 
every {\it smooth} positive definite $Z^{\pi_1}$-manifold $M$ 
is  TOP-split, because for every $m$ the $Z^{\pi_1}$-manifold $M^{(m)}$ is a  positive definite smooth 4-manifold 
 (see Lemma~2.1 later)  and hence the intersection form of $M^{(m)}$ is represented by a block sum of copies 
of $(1)$ by Donaldson's theorem in \cite{D}, where note that the intersection form of $M^{(m)}$ is 
identical to the intersection form of 
a closed simply connected smooth 4-manifold obtained from  $M^{(m)}$ by killing $\pi_1(M^{(m)})=Z$.  

Thus,  every {\it smooth} $Z^{\pi_1}$-manifold $M$ is TOP-split, 
because every indefinite $Z^{\pi_1}$-manifold is seen in \cite{KA} to be TOP-split.

As another consequence (shown in \cite{KA2}), every {\it smooth} $S^2$-knot in every 
{\it smooth} closed 
simply connected 4-manifold is topologically unknotted if the knot group is an infinite cyclic group. 

For a positive definite $Z^{\pi_1}$-manifold $M$, assume that there are 
$\Lambda$-generators 
$\tilde x_i\, (i=1,2,\dots,n)$ of $H_2(\widetilde M;Z)$ 
with $||\tilde x_i||^2 \leqq 2$ 
for all $i$.   Let $m$ be  an integer such that  
\[m\geqq \max\{e({\tilde x_1})+1,e({\tilde x_2})+1,
\dots,e({\tilde x_n})+1,\lambda\}\] 
for some winding degree $\lambda$ on $M$. Then the elements 
$\tilde x_i^{(m)}\in H_2(M^{(m)};Z)\, (i=1,2,\dots,n)$ induced from $\tilde x_i\, (i=1,2,\dots,n)$  by the 
covering projection homomorphism $H_2(\widetilde M;Z)\to H_2(M^{(m)};Z)$  
form $\Lambda^{(m)}$-generators of $H_2(M^{(m)};Z)$ with  identical square length  
$||\tilde x_i^{(m)}||^2=||\tilde x_i||^2$
for all $i$, so that there are $Z$-generators 
$x'_{i'}\,(i'=1,2,\dots, n')$ for $H_2(M^{(m)};Z)$ such that 
$||x'_{i'}||^2\leqq 2$  
for all $i'$, meaning that $M$ is TOP-split by Corollary~1.2. 
Thus, the following corollary is obtained from Theorem~1.1 (4) and  
the observation just before Corollary~1.2 (see Corollary~3.3 later).  

\phantom{x}

\noindent{\bf Corollary~1.3.} A positive 
definite $Z^{\pi_1}$-manifold $M$ is TOP-split if 
there are $\Lambda$-generators 
$\tilde x_i\, (i=1,2,\dots,n)$ of $H_2(\widetilde M;Z)$  
such that $||\tilde x_i||^2 \leqq 2$ for all $i$. 

\phantom{x}

The following theorem gives another criterion that a positive definite 
$Z^{\pi_1}$-manifold with  standard $Z$-intersection form is TOP-split.

\phantom{x}

\noindent{\bf Theorem~1.4.} Let $M$ be a positive 
definite $Z^{\pi_1}$-manifold such that there is a 
$Z$-basis $e_i \, (i=1,2,\dots,n)$ of $H_2(M;Z)$ 
with $\mbox{Int}^M(e_i,e_j)=\delta_{ij}$ for all 
$i, j$. Then the following conditions (0), (1), (2) and (3) on $M$ are mutually 
equivalent:

\medskip

\noindent(0) The $Z^{\pi_1}$-manifold $M$ is TOP-split. 

\noindent(1) For a $\Lambda$-basis $\tilde x_i\, (i=1,2,\dots,n)$ 
of $H_2(\widetilde M;Z)$,  there are elements $a_{ij}(t)\in \Lambda\, 
(i,j=1,2,\dots,n)$ with $\Lambda$-intersection number 
\[\mbox{Int}^{\widetilde M}_{\Lambda}(\tilde x_i, \tilde x_j)
=\sum_{k=1}^n a_{ik}(t^{-1})a_{jk}(t)\]
for every $i$ and $j$.

\noindent(2) Every minimal element $\tilde x\in H_2(\widetilde M;Z)$ has  square length $||\tilde x||^2=1$.  

\noindent(3) The following conditions (3${}_1$) and (3${}_2$) are satisfied: 

\noindent(3${}_1$) For every  element $\tilde x \in 
H_2(\widetilde M;Z)$, there are elements $a_i(t)\in \Lambda\, 
(i=1,2,\dots,s)$ such that the $\Lambda$-square length 
\[||\tilde x||^2_{\Lambda}= \sum_{i=1}^s a_i(t)a_i(t^{-1}).\]

\noindent(3${}_2$)  There are mutually distinct 
(up to multiplications of 
the units of $\Lambda$) elements $\tilde y_i\, (i=1,2,\dots,n-1)$ in 
$H_2(\widetilde M;Z)$ such that the square length  
$||\tilde y_i||^2=1$ for all $i$.

\phantom{x}

For example, every $Z^{\pi_1}$-manifold $M$ with  second 
Betti number $\beta_2(M)=1$ can be converted to  a positive definite $Z^{\pi_1}$-manifold 
with standard $Z$-intersection form by changing the orientation if necessary. Such a manifold  
is TOP-split by Theorem~1.4 since the $\Lambda$-intersection matrix  
on $H_2(\widetilde M;Z)\cong \Lambda$ is $(\pm1)$. It appears to be unknown 
{\it whether every positive definite} $Z^{\pi_1}$-{\it manifold} $M$ {\it with} 
$\beta_2(M)=2$ {\it or} $3$ {\it is TOP-split}, where we note that 
every positive definite $Z$-form of rank up to $7$ is  known to be 
standard (see \cite{MH}). 

We come back to 
the  Hambleton-Teichner-Freedman-Quinn $Z^{\pi_1}$-manifold $M_L$.
The matrix
\[L=\left(
\begin{array}{cccc}
1+f+f^2 & f+f^2 & 1+f& f\\
f+f^2 & 1+f+f^2 & f & 1+f\\
1+f & f & 2 & 0\\
f & 1+f & 0 & 2
\end{array}
\right)\]
with $f=t+t^{-1}$ is the Hambleton-Teichner matrix given in \cite{HT}.  
Since the size of the matrix $L$ is $4$, the $Z$-intersection form 
$\mbox{Int}^{M_L}$ on $H_2(M_L;Z)$ is the standard form. 
Let $\tilde x_i\, (i=1,2,3,4)$ be the 
$\Lambda$-basis of $H_2(\widetilde M_L;Z)$ giving 
$L$ as the $\Lambda$-intersection matrix. 
It is shown by Theorem~1.4 (3${}_1$) that $M_L$ is not TOP-split  
since 
$||\tilde x_i||^2=1+f+f^2\, (i=1,2)$ 
cannot be written as a sum $\sum_{i=1}^s a_i(t)a_i(t^{-1})$.  
The $Z^{\pi_1}$-manifold $M_L$ has further properties 
which are given below:

\phantom{x}

\noindent{\bf Theorem~1.5.}  For the Hambleton-Teichner matrix $L$, 
we have the following properties (1) and (2).  

\medskip

\noindent(1) The $\Lambda$-basis $\tilde x_i,\, (i=1,2,3,4)$ of 
$H_2(\widetilde M_L;Z)$ are minimal 
elements such that the square length 
$||\tilde x_i||^2$ is $3$ for $i=1,2$ and $2$ for $i=3,4$.

\noindent(2) Let $V_G$ be the $Z$-intersection matrix on any $Z$-free 
subgroup $G$ of $H_2(\widetilde M_L;Z)$ of finite rank. Then 
the determinant $\mbox{det} V_G$ of $V_G$ is greater than $1$.  

\medskip

In particular, we see that the criteria of Theorem~1.4 are not satisfied for $M_L$. 

\phantom{x}

The following corollary produces an infinite family of 
orthogonally indecomposable 
unimodular odd symmetric $Z$-forms as a corollary to Theorem~1.5.  

\phantom{x}

\noindent{\bf Corollary~1.6.}  For the Hambleton-Teichner matrix $L$ 
and every integer $m\geqq 3$, 
we have the following properties (1) and (2). 

\medskip

\noindent(1) The $\Lambda^{(m)}$-basis 
$\tilde x_i^{(m)},\, (i=1,2,3,4)$  of 
$H_2(M_L^{(m)};Z)$ are minimal 
elements such that the square length  
$||\tilde x_i^{(m)}||^2$ is $3$ for 
$i=1,2$ and $2$ for $i=3,4$.

\noindent(2) The $Z$-intersection form on $H_2(M_L^{(m)};Z)\cong Z^{4m}$ 
is an orthogonally indecomposable unimodular odd definite symmetric 
$Z$-form.

\phantom{x}

It was shown in \cite{FHMT} by  a different method that the 
$Z$-intersection form 
on $H_2(M_L^{(m)};Z)$ in (2) is not standard for every $m\geqq 3$ and 
orthogonally indecomposable for $m=3,4$. The proof of Corollary~1.6 (2) will be done for $m\geqq 5$. 
The form  in (2) is standard for 
$m\leqq 2$ and isomorphic to $\Gamma_{12}$ for $m=3$ by the classification of 
orthogonally indecomposable unimodular definite symmetric $Z$-forms of rank $\leqq 16$ (see \cite{MH}). 
On the other hand, the form in (2) 
must be different from $\Gamma_{4m}$ for every $m\geqq 4$, 
because the square 
length of every minimal element of $\Gamma_{4m}$ is $2$ or $m$ (see \cite{MH}). 

As another note on Corollary~1.6,  the elements $\tilde x^{(1)}_i\in H_2(M_L;Z)$ and 
$\tilde x^{(2)}_i\in H_2(M^{(2)};Z)$ for $i=1,2$ have the equalities
\[||\tilde x^{(2)}_i||^2=||\tilde x^{(1)}_i||^2-2=3,\]
which give a concrete example that the inequality of Theorem~1.1 (3) 
does not hold in general unless $M$ is TOP-split.

Given a Hermitian $\Lambda$-matrix $A$ with $\mbox{det}A=1$, we can construct a {\it smooth} 
compact connected oriented  4-manifold $E$ from  the $4$-dimensional solid torus $S^1\times D^3$ 
by attaching some $2$-handles to the boundary $\partial (S^1\times D^3)=S^1\times S^2$  
such that $\pi_1(E;Z)\cong Z$ and the matrix $A$ is a $\Lambda$-intersection matrix on 
$H_2(\widetilde E;Z)$, where $\widetilde E$ denotes the infinite cyclic connected covering space of $E$ (see \cite{FQ}).  
The boundary $\partial E=B$ of $E$ has the same homology as $S^1\times S^2$ and hence is called 
a {\it homology handle}. Since $\mbox{det}A=1$, the boundary $\widetilde B= \partial \widetilde E$ of  
$\widetilde E$ has the trivial homology $H_1(\widetilde B;Z)=0$. 
In case $A$ is the Hambleton-Teichner matrix $L$, it is observed in 
\cite{FHMT} that the homology handle $B$ cannot bound a {\it smooth} rational 
homology circle $W$ by using \cite{D2} instead of \cite{D}, 
which is generalized as follows:\footnote{It is assumed in \cite{FHMT} that 
the natural homomorphism: 
$H_1(B;Z)\to H_1(W;Z)/(\mbox{torsion})$ is an isomorphism, but this restriction can 
be removed.}

\phantom{x}

\noindent{\bf Observation~1.7.} For the Hambleton-Teichner matrix $L$,  
the disjoint union $n B$ of $n$ copies of the homology handle $B$ 
for any $n\geqq 1$ cannot bound a smooth 
compact oriented 4-manifold $W$ with $H_2(W;Q)=0$. 

\phantom{x}

It is unknown {\it whether 
the homology handle \mbox{$B$} can  bound a smooth compact oriented 
\mbox{4}-manifold \mbox{$W$} with an infinite cyclic covering space 
\mbox{$\widetilde W$} 
such that} 
\[\dim_Q H_2(\widetilde W;Q)<+\infty,\] in other words,  
{\it whether the homology handle \mbox{$B$} represents a trivial  
element of the \mbox{$\widetilde H$}-cobordism group} 
$\widetilde\Omega(S^1\times S^2)$ in  \cite{Kold0}.  
It is also observed in \cite{FHMT} that some $B$ can be taken as 
the Dehn surgery manifold  with coefficient $0$ along a knot $K$ in the 3-sphere $S^3$ with trivial Alexander 
polynomial. With an idea of Cochran-Lickorish \cite{CL}, we obtain:

\phantom{x}

\noindent{\bf Observation~1.8.} For the Hambleton-Teichner matrix $L$, 
the knot $K$ cannot be converted into the unknot 
by changing positive crossings into negative crossings. 

\phantom{x}

In \S~2, several preliminaries on a winding degree are provided. 
In \S~3, the proofs of Theorems~1.1 and 1.4 are given.   
In \S~4, the proofs of  Theorem~1.5, Corollary~1.6 and  Observations~1.7 and 1.8 are given. 

\phantom{x}

\noindent{\bf 2. Several preliminaries on the winding degree}

For a $Z^{\pi_1}$-manifold $M$,  
let $\tilde x\in H_2(\widetilde M;Z)$ be a non-zero element.  
We assert that the square length $||\tilde x||^2>0$. To see this, 
the following lemma is proved  (though it was the fact used in  \cite{KA2}). 

\phantom{x}

\noindent{\bf Lemma~2.1}
If a $Z^{\pi_1}$-manifold $M$ is  
positive definite, then $M^{(m)}$ is also  positive definite 
for any $m$. 

\phantom{x}

\noindent{\bf Proof.} There  are two ways to see that $\beta_2(M^{(m)})=m\beta_2(M)$.  
One way is to use the Euler characteristic identity $\chi(M^{(m)})=m\chi(M)$.  By the Betti numbers 
 $\beta_d(M^{(m)})=\beta_d(M)=1$ for $d=0,1$, we have $\beta_2(M^{(m)})=m\beta_2(M)$. 
The other way is to use that $H_2(\widetilde M^{(m)};Z)$ is a $\Lambda^{(m)}$-module of 
rank $\beta_2(M)$, showing that $\beta_2(M^{(m)})=m\beta_2(M)$.   
If $\sigma(M)=\beta_2(M)$, then  the signature identity $\sigma(M^{(m)})=m\sigma(M)$ 
shows that $\sigma(M^{(m)})=\beta_2(M^{(m)})$.  $\square$

\phantom{x}

Note that the image $\tilde x^{(m)}\in H_2(M^{(m)};Z)$ of $\tilde x$ 
under the covering projection homomorphism 
$H_2(\widetilde M;Z)\to H_2(M^{(m)};Z)$ is not zero for a large $m$. 
Hence, we have the square length $||\tilde x^{(m)}||^2>0$ for a large $m$. 
The $\Lambda^{(m)}$-intersection number 
$\mbox{Int}^{M^{(m)}}_{\Lambda^{(m)}}(\tilde x^{(m)}, \tilde y^{(m)})$ 
of the elements 
$\tilde x^{(m)}, \tilde y^{(m)} \in H_2(M^{(m)};Z)$ is calculated 
as follows:
\[\mbox{Int}^{M^{(m)}}_{\Lambda^{(m)}}(\tilde x^{(m)}, \tilde y^{(m)})
=\sum_{s=0}^{m-1}\left(\sum_{i=-\infty}^{+\infty}
\mbox{Int}^{\widetilde M}(t^{s+mi}\tilde x, \tilde y)t^s\right)\in \Lambda^{(m)}.\] 
In particular, the $Z$-intersection number 
$\mbox{Int}^{M^{(m)}}(\tilde x^{(m)}, \tilde y^{(m)})$ is given as follows:
\[\mbox{Int}^{M^{(m)}}(\tilde x^{(m)}, \tilde y^{(m)})
=\sum_{i=-\infty}^{+\infty}
\mbox{Int}^{\widetilde M}(t^{mi}\tilde x, \tilde y)\in Z.\]
By definition, for the exponent $e(\tilde x)$ of $\tilde x$,  which is 
the highest degree of the $\Lambda$-square length
$||\tilde x||^2_{\Lambda}=\mbox{Int}^{\widetilde M}_{\Lambda}(\tilde x,\tilde x)$, 
it  is seen that the square length $||\tilde x||^2$ 
is identical to the square length $||\tilde x^{(m)}||^2$ for any integer 
$m\geqq e(\tilde x)+1$, so that  if $\tilde x\ne 0$ and 
$m\geqq e(\tilde x)+1$, then we have 
$||\tilde x||^2=||\tilde x^{(m)}||^2>0$. This shows that the square length 
 $||\tilde x||^2>0$ as asserted. 

\phantom{x}

\noindent{\bf Definition.}
A 3-{\it sphere leaf} of a TOP-split $Z^{\pi_1}$-manifold $X$ is 
a 3-sphere submanifold $V$ of $X$ corresponding to the 3-sphere $1\times S^3$ 
under a homeomorphism $X\cong S^1\times S^3\# X_1$ for some simply connected 
4-manifold $X_1$. 

\phantom{x}

Let $P$ be a 2-sphere (embedded) in $X$. 
We assume that the intersection $L=P\cap V$ is a closed oriented 
possibly disconnected 1-manifold unless it is empty. 
Let $D_i\, (i=1,2,\dots,r)$ be the connected regions of $P$ divided by $L$. 
Let $p_i$ be a fixed interior point of $D_i$. 
Let $\alpha_{ij}$ be an oriented arc in $P$ joining the point $p_i$ 
to the point $p_j$ transversely meeting $L$. The absolute value 
$|\mbox{Int}^P(\alpha_{ij},L)|$ of the $Z$-intersection number 
$\mbox{Int}^P(\alpha_{ij},L)$ is independent of any choices of  the points $p_i\in D_i$, $p_j\in D_j$ 
and the oriented arc $\alpha_{ij}$. 
The maximal number $|\mbox{Int}^P(\alpha_{ij},L)|$ for all $i,j$ is determined  only by the 2-sphere $P$ in $X$ 
and the 3-sphere leaf $V$ and called the 
{\it winding index} of the 2-sphere $P$ in $X$ with respect to 
the 3-sphere leaf $V$ and denoted by $\delta=\delta(P;V,X)$.  A  meaning of the winding index $\delta$ is 
as follows. 
Let $X_V$ be the fundamental region  of the infinite cyclic covering space $\widetilde X$, 
obtained from $X$ by splitting along $V$.  Note that $\widetilde X$ is the union of  the $t^i$-shifts $t^i(X_V)\, 
(i=0,\pm1,\pm2,\dots)$ of $X_V$. 
Then the winding index $\delta$ is the maximal number of the interiors of $t^i(X_V)\,(i=0,\pm1,\pm2,\dots)$
meeting a fixed lifting 2-sphere  $\widetilde P$ of  $P$ to $\widetilde X$. 
This notion was defined in \cite{KA2}.  

A homological version of the winding index for a 2-sphere $P$  in   a TOP-split $Z^{\pi_1}$-manifold $X$ with 
a 3-sphere leaf $V$  as follows:

The homology $H_2(\widetilde X;Z)$ is 
a free $\Lambda$-module with a $\Lambda$-basis $v_i\, (i=1,2,\dots,n)$ 
represented by 2-cycles in $X_V$. 
For an element $\tilde x\in  H_2(\widetilde X;Z)$, 
the $\Lambda$-intersection number $\mbox{Int}_{\Lambda}(\tilde x,v_i)$ 
is an element of $\Lambda$, i.e., an integral  Laurent polynomial in $t$. 
The maximal and minimal 
degrees  of $\mbox{Int}_{\Lambda}(\tilde x,v_i)$ in $t$ for all $i$ are independent of a choice 
of  such a $\Lambda$-basis $v_i\, (i=1,2,\dots,n)$ and  denoted by 
$\max \mbox{deg}(\tilde x;V,X)$ and $\min \mbox{deg}(\tilde x;V,X)$, 
respectively. The difference 
$\delta^h(\tilde x;V,X)=\max \mbox{deg} (\tilde x;V,X)-\min \mbox{deg} (\tilde x;V,X)$ is independent of 
any covering translations of $\tilde x$ and called the 
{\it homological winding index} of  $\tilde x$ in $X$ with respect to the 
3-sphere leaf $V$. The {\it homological winding index}$\delta^h(P;V,X)$ of the 2-sphere $P$ in $X$ 
with respect to the 3-sphere leaf $V$ is the number 
 $\delta^h([\widetilde P];V,X)$ for the homology class $[\widetilde P]\in H_2(\widetilde X;Z)$ of 
 any lifting 2-sphere $\widetilde P$ of $P$ to $\widetilde X$.  

It is direct to see that 
\[\delta^h(P;V,X)\leqq \delta(P;V,X).\]

The notion of a winding degree on a positive definite $Z^{\pi_1}$-manifold 
$M$ is defined  from the notion of a homological winding index as follows:
  
For an integer $u\geqq 1$, 
let $Q(u)=\#_u{{\mathbf C}P}^2$ 
and $\overline Q(u)=\#_u{\overline{{\mathbf C}P}}^2$ 
be the $u$-fold connected sum of 
the complex projective planes $Q={{\mathbf C}P}^2$ and 
$\overline Q={\overline{{\mathbf C}P}}^2$ with signatures $+1$ and $-1$, 
respectively. 
Let  $uP$ be the disjoint union of 
the 2-spheres $P_k=\overline{{\mathbf C} P}^1_k\, 
(k=1,2,\dots, u)$ in the connected summands 
$\overline Q_k={\overline{{\mathbf C}P}}^2_k\, (k=1,2,\dots, u)$ of $\overline Q(u)$. 
A tubular neighborhood $N\subset X$ of a generator circle of $\pi_1(X)\cong Z$ 
is called a {\it solid tube generator} of $X$ which is unique up to ambient isotopies of $X$.  
A {\it circle union} $X'\circ X''$ of 
two $Z^{\pi_1}$-manifolds $X',X''$ is a $Z^{\pi_1}$-manifold obtained from the 
exteriors 
$\mbox{cl}(X'\setminus N')$ and $\mbox{cl}(X''\setminus N'')$ for solid 
tube generators $N'$ and $N''$ of $X'$ and $X''$, respectively, by identifying 
the boundary $\partial N'$ with  boundary $\partial N''$ by 
an orientation-reversing homeomorphism. 
For a positive definite $Z^{\pi_1}$-manifold $M$, we consider 
a $Z^{\pi_1}$-manifold $X=\overline Q(u)\# Q(v)\# M$ 
for integers  $u\geqq 1$ and $v\geqq 0$. Since $X$ is an indefinite 
$Z^{\pi_1}$-manifold,  $X$ is TOP-split by \cite{KA}. 
Consider $X$ as a circle union $X'\circ X''$ 
with following conditions (i)-(iii).

\medskip 

\noindent(i)\, The $Z^{\pi_1}$-manifold $X'$ is a TOP-split 
$Z^{\pi_1}$-manifold 
with  a 3-sphere leaf $V'$ with $V'\cap \mbox{cl}(X'\setminus N')=B'$ a 3-disk. 

\noindent(ii)\, For a solid tube generator $N\subset M$, 
we have  an inclusion 
$\mbox{cl}(M\setminus N)\subset \mbox{cl}(X'\setminus N')$ inducing an 
isomorphism on the infinite cyclic fundamental groups.  

\noindent(iii)\, The $Z^{\pi_1}$-manifold
$X''$ is a positive definite $Z^{\pi_1}$-manifold. 

\medskip

This circle union splitting $X'\circ X''$ of $X$ always exists although it is not unique.  
For example, let $X'=\overline Q(u)\# Q(v')\# M$ (which is TOP-split by \cite{KA}) and $X''=Q(v'')\#S^1\times S^3$ 
for any sum $v=v'+v''$ which give  a desired circle union $X'\circ X''$  of $X$ for any 3-sphere leaf $V'$ of $X'$. 
In fact,  for a solid  tube generator $N$ of $M$, there is a 3-sphere leaf 
$V'$ of $X'$ with $N\cap V'$ a 3-ball.  Let  $N'= \mbox{cl}(N\setminus c(\partial N))$ for a boundary collar  
$c(\partial N)$ of $\partial N$ in $N$. Then $X'$ is a union of 
$\mbox{cl}(M\setminus N)\cup (c(\partial N)\# \overline Q(u)\#Q(v'))\cup N'$. 
In this decomposition, the conditions  (i) and (ii) are satisfied. The condition (iii) is clearly satisfied.  

Then note that there is a canonical isomorphism 
\[H_2(\widetilde X;Z)\cong H_2(\widetilde X';Z)\oplus H_2(\widetilde X'';Z).\] 
For a connected lift $\widetilde P_k$ of the 2-sphere $P_k$ to $\widetilde X$, let $[\widetilde P_k]'$ 
be the projection image of the homology class 
$[\widetilde P_k]\in H_2(\widetilde X;Z)$ into the direct summand $H_2(\widetilde X';Z)$. 
The maximum of the homological winding index $\delta^h([\widetilde P_k]';V',X')$ for all 
$k\,(k=1,2,\dots, u)$ is denoted by $\delta^h([uP]';V',X')$. 

\phantom{x}

\noindent{\bf Definition.}
A {\it winding degree} on a positive definite $Z^{\pi_1}$-manifold 
$M$, denoted by $\lambda$,  is the non-negative integer 
$\delta^h([uP]';V',X')$ given by a choice of  integers $u\geqq 1$ and $v\geqq 0$ for the connected sum 
$X=\overline Q(u)\# Q(v)\# M$,  a choice of any circle union splitting $X'\circ X''$ of $X$ with  properties 
(i)-(iii), and a choice of a 3-sphere leaf $V'$ of $X'$. 

\phantom{x}

A  winding degree $\lambda=\delta^h([uP]';V',X')$ is an amount that measures the possible range.
of the homology classes $[\widetilde P_k]\in H_2(\widetilde X;Z)$ in the decompositions of $\widetilde X'$ into 
the $t$-power shifts $t^k(M_{V'})\, (i=\pm 1,\pm2,\dots)$ of the fundamental region $X'_{V'}$. 
Let $\lambda_{\mbox{min}}$ be the minimum of all winding degrees $\lambda$ on a positive definite 
$Z^{\pi_1}$-manifold$M$ which is, by definition, a topological invariant of $M$. 

Let $P^E_k=P_k\cap \mbox{cl}(X'\setminus N')$ be a proper surface in the compact manifold 
$\mbox{cl}(X'\setminus N')$ whose boundary is an oriented link $L^P_k$ with orientation induced from $P^E_k$ such that 
$L^P_k$ is in a 3-ball $B''\subset \partial N'\setminus \partial B'$ by an ambient isotopic deformation of  $P^E_k$. 
Let $D_k$ be a connected Seifert surface for $L^P_k$ in $B''$.  
Let $PD_k=P^E_k\cup (-D_k)$ be the closed connected oriented surface in $X'$, and 
$\widetilde PD_k$ a connected lift of $PD_k$ to $\tilde X'$. 
Note that  the homology class $[\widetilde P_k]'\in H_2(\widetilde X';Z)$ 
which is the projection image of the homology class 
$[\widetilde P_k]\in H_2(\widetilde X;Z)$ into  the direct summand $H_2(\widetilde X';Z)$  is written as 
\[[\widetilde P_k]'=[\widetilde PD_k]=t^r(c_0 +c_1 t+\dots +c_d t^d),\]
where $r=\min \mbox{deg}([\widetilde P_k]';V',X')$,  
$r+d=\max \mbox{deg}([\widetilde P_k]';V',X')$, and 
$c_i\,(i=0,1,\dots,d)$ 
are homology classes in $H_2(\widetilde X';Z)$ represented by 2-cycles in the 
fundamental region $X'_{V'}$ with $c_0\ne 0$ and $c_d\ne 0$. 

Let $\widetilde S^3$ be a fixed 3-sphere lift of $V'$ to 
$\widetilde X'$.  For any integer $j$, let 
$\widetilde L_{k,j}=t^{-j}(\widetilde PD_k\cap t^j \widetilde S^3)$ 
be an oriented link (with orientation determined by 
the orientations of $\widetilde PD_k$, $\widetilde S^3$ and $\widetilde X$) 
in $\widetilde S^3$ unless it is empty. 
The following lemma is used in our argument.  

\phantom{x}

\noindent{\bf Lemma~2.2.}  
Let $F$ be a closed oriented surface in $X'\setminus N'$ 
with $F\cap PD_k=\emptyset$, and $\widetilde F$ the preimage of $F$ under the projection 
$\widetilde X'\to X'$.  
Assume that the surface $\widetilde F$ meets the 3-sphere $\widetilde S^3$ 
as a knot $\widetilde K$.  Then for any integer $j$ with 
$j\leqq \min \mbox{deg}([\widetilde P_k]';V',X')$ or 
$j>\max \mbox{deg}([\widetilde P_k]';V',X')$ with $\widetilde L_{k,j} \ne\emptyset$,  
the linking number ${\mbox{Link}}^{\widetilde S^3}(\widetilde K,\widetilde L_{k,j})$ in 
the 3-sphere $\widetilde S^3$ is $0$. 

\phantom{x}

\noindent{\bf Proof.} 
For any $j$ with $\widetilde L_{k,j}\ne\emptyset$, we construct a closed oriented 
surface $t^j(C'_{k,j})\cup (-C''_{k,j})$ in $\widetilde X'$ where 
$t^j(C'_{k,j})$ is  a Seifert surface  of the link 
$t^j (\widetilde L_{k,j})$ in the 3-sphere $t^j (\widetilde S^3)$ and 
$C''_{k,j}$ is a compact surface in $\widetilde PD_k$ bounded by 
$t^j(\widetilde L_{k,j})$. 
If $j\leqq r$ or $j>r+d$, then it is possible 
to choose $C''_{k,j}$ so that the surface 
$t^j(C'_{k,j})\cup (-C''_{k,j})$ is null-homologous in $\widetilde X'$. 
By using the fundamental region $X'_{V'}$, choose a compact 4-submanifold 
$X'_{J}=\cup_{i=-J}^{J} t^i(X'_{V'})$ of $\widetilde X'$
for a sufficiently large integer $J$ to contain the surface 
$t^j(C'_{k,j})\cup (-C''_{k,j})$ and the 3-sphere $t^j (\widetilde S^3)$ 
in the interior. 
Let $\hat F^c$ be a closed oriented surface 
obtained from the compact surface $F^c=\widetilde F\cap X_J$ by 
adding  surfaces in $\partial X_J$ which are translations of a Seifert surface of 
$\widetilde K$ in $\widetilde S^3$. 
Then the $Z$-intersection number 
$\mbox{Int}^{\widetilde X}(\hat F^c, t^j(C'_{k,j})\cup (-C''_{k,j}))$  is zero. 
Since $C''_{k,j}\cap \hat F^c=\emptyset$, 
the $Z$-intersection number 
$\mbox{Int}^{t^j(\widetilde S^3)}(t^j(\widetilde K), t^j(C'_{k,j}))$ 
is zero and hence ${\mbox{Link}}^{\widetilde S^3}(\widetilde K,\widetilde L_{k,j})=0$. $\square$

\phantom{x}

Throughout the remainder of this section, an estimate of a winding degree from 
a $\Lambda$-intersection matrix for an odd positive definite $Z^{\pi_1}$-manifold 
is explained. 
Let $A=(a_{ij}(t))$ be a $\Lambda$-intersection matrix of size $n$ 
on an odd positive definite $Z^{\pi_1}$-manifold $M$. 
It is noted that $a_{ij}(t)=a_{ji}(t^{-1})$ for all $i,j$. 
Let $\tilde x_i\,(i=1,2,\dots,n)$ be the  $\Lambda$-basis for 
$H_2(\widetilde M;Z)$ 
giving the matrix $A$, and $\tilde x'_i\,(i=1,2,\dots,n)$ the dual 
$\Lambda$-basis, i.e., the $\Lambda$-basis for $H_2(\widetilde M;Z)$ with 
$\mbox{Int}^{\widetilde M}_{\Lambda}(\tilde x_i,\tilde x'_j)=
\mbox{Int}^{\widetilde M}_{\Lambda}(\tilde x'_j,\tilde x_i)=\delta_{ij}$ 
for all $i, j$, whose $\Lambda$-intersection matrix is given by the inverse 
matrix $A^{-1}=(b_{ij}(t))$. The following identities 
\begin{eqnarray*}
\tilde x_j&=&\sum_{k=1}^n a_{kj}(t) \tilde x'_k\quad (j=1,2,\dots,n),\\
\tilde x'_j&=&\sum_{k=1}^n b_{kj}(t) \tilde x_k\quad (j=1,2,\dots,n)
\end{eqnarray*}
are easily established. Consider the following unique splittings of the 
Laurent polynomials $a_{ii}(t)$ and $b_{ii}(t)$ in $t$:
\[a_{ii}(t)=\varepsilon^a_i+a'_{ii}(t)+a'_{ii}(t^{-1}),\quad 
b_{ii}(t)=\varepsilon^b_i+b'_{ii}(t)+b'_{ii}(t^{-1})\] 
where $\varepsilon^a_i$ and $\varepsilon^b_i$ are taken $0$ or $1$ and 
$a'_{ii}(t)$ and $b'_{ii}(t)$ are elements in $\Lambda$ with  
non-negative constant terms and without any negative powers of $t$. 
Let 
\[
{\tilde a_{ij}(t)=
\left\{
\begin{array}{ll}
a_{ij}(t) &\quad i\ne j\\
a'_{ii}(t) &\quad i=j,
\end{array}
\right.}
\qquad
{\tilde b_{ij}(t)=
\left\{
\begin{array}{ll}
b_{ij}(t) &\quad i\ne j\\
b'_{ii}(t) &\quad i=j.
\end{array}
\right.}
\]
Further, for a double index element 
$f_{i j}(t)\in \Lambda\, (i,j=1,2,\dots,n)$, let 
\[
f_{i j}^0(t)=
\left\{
\begin{array}{ll}
0 &\quad i>j\\
f_{i j}(t) &\quad i\leqq j.
\end{array}
\right.
\qquad
{f_{i j}^{0*}(t)=
\left\{
\begin{array}{ll}
0 &\quad i<j\\
f_{i j}(t) &\quad i\geqq j.
\end{array}
\right.}
\]

Let $\max\lambda(A)$ and $\min\lambda(A)$ be respectively the maximal degree and 
the minimal degree of 
the following Laurent polynomials in $t$:
\[1, \quad \ a_{ij}(t)\, (i<j), \quad b_{ij}(t)\, (i<j),\quad a'_{ii}(t^{-1})
,\quad b'_{ii}(t),\quad 
c_{ij}(t)=\sum_{k=1}^{\min\{i,j\}}\tilde a_{ik}(t)\cdot\tilde b_{kj}(t)\] 
for all $i,j$.  Let $\lambda(A)=\max\lambda(A)-\min\lambda(A)$. 
Then an estimate of a winding degree is done as follows:

\phantom{x}

\noindent{\bf Lemma~2.3.} For an odd positive definite 
$Z^{\pi_1}$-manifold $M$ with $\Lambda$-intersection matrix $A=(a_{ij}(t))$, 
there is a 
winding degree $\lambda$ on $M$ such that $\lambda\leqq \lambda(A)$. 

\phantom{x}

\noindent{\bf Proof.} 
For the orientation-reversed manifold $-M$ of $M$, it is noted that 
any circle union $Y=M\circ(-M)$ is TOP-split as 
$S^1\times S^3\#Q(n)\#\overline Q(n)$ 
because odd indefinite forms are diagonal and the connected sum $Y_1=M_1\#(-M_1)$ 
for the simply connected 
manifold $Y_1$ obtained from  $Y$ by a 2-handle surgery killing 
$\pi_1(Y)\cong Z$ is homeomorphic to $Q(n)\#\overline Q(n)$ by using 
the vanishing of Kirby-Siebenmann obstruction (see \cite{FQ}) and 
the manifold $Y$ is a TOP-split $Z^{\pi_1}$-manifold as it is discussed from now. 
For the $\Lambda$-basis $\tilde x_i\,(i=1,2,\dots,n)$ and its dual $\Lambda$-basis $\tilde x'_i\,(i=1,2,\dots,n)$ 
for $H_2(\widetilde M;Z)$ representing $A$ and $A^{-1}$, respectively, 
let ${\tilde x}^-$ and $\tilde x^{\prime-}_i$ correspond to $\tilde x$ and 
$\tilde x'_i$ in the direct summand $H_2(-\widetilde M;Z)$ of $H_2(\widetilde Y;Z)$ for every $i$, respectively. 
Let $\tilde y_i={\tilde x'}_i+\tilde x^{\prime-}_i\in 
H_2(\widetilde Y;Z)\, 
(i=1,2,\dots, n)$.  Note that 
$\mbox{Int}^{\widetilde Y}_{\Lambda}(\tilde y_i,\tilde y_j)=0$ 
and 
$\mbox{Int}^{\widetilde Y}_{\Lambda}(\tilde y_i,\tilde x^-_j)=-\delta_{ij}$ 
for all $i,j$. The elements $\tilde z_i\, (i=1,2,\dots,n)$ in 
$H_2(\widetilde Y;Z)$ 
are constructed as follows:
\begin{eqnarray*}
\tilde z_1 &=& \tilde x^-_1- a'_{11}(t^{-1})\tilde y_1,\\
\tilde z_2 &=& \tilde x^-_2- a_{21}(t^{-1})\tilde y_1-a'_{22}(t^{-1})\tilde y_2,\\
           &\dots&\\
\tilde z_n &=& \tilde x^-_n- a_{n1}(t^{-1})\tilde y_1-a_{n2}(t^{-1})\tilde y_2-
\dots-a'_{nn}(t^{-1})\tilde y_n.
\end{eqnarray*}
It is checked that 
\[\mbox{Int}^{\widetilde Y}_{\Lambda}(\tilde y_i,\tilde z_j)=-\delta_{ij}, \quad 
\mbox{Int}^{\widetilde Y}_{\Lambda}(\tilde z_i,\tilde z_j)=-\varepsilon^a_i \delta_{ij}\]
for all $i,j$, so that 
$\tilde y_i,\tilde z_j\, (i,j=1,2,\dots,n)$ form a $\Lambda$-basis 
for $H_2(\widetilde Y;Z)$ with integral $\Lambda$-intersection matrix 
and thus,  the $Z^{\pi_1}$-manifold $Y$ is TOP-split.  
Since $Y$ is odd and has the trivial Kirby-Siebenmann obstruction, there is an 
identification 
\[Y=S^1\times S^3\#Q(u)\#\overline Q(u)\] 
where every component of the 2-spheres $nP$ in $\overline Q(n)$ represents a $Z$-linear 
combination 
of the $\Lambda$-basis $\tilde y_i,\tilde z_j\, (i,j=1,2,\dots,n)$ of 
$H_2(\widetilde Y;Z)$. 
Let $X=Y\circ M=M\circ(-M)\circ M=M\circ Y^*$ with $Y^*=(-M)\circ M$.
A $\Lambda$-basis $\tilde y^*_i,\tilde z^*_j\, (i,j=1,2,\dots,n)$ 
for $H_2(\widetilde Y^*;Z)$ is given as follows:

Let 
$\tilde y^*_i=\tilde x^-_i+\tilde x^*_i$ where $\tilde x^-_i$ and 
$\tilde x^*_i$
correspond to $\tilde x^-_i\in H_2(-\widetilde M;Z)$ 
and $\tilde x_i\in H_2(\widetilde M;Z)$, 
respectively. Then we have 
$\mbox{Int}^{\widetilde Y^*}_{\Lambda}(\tilde y^*_i,\tilde y^*_j)=0$ 
for all $i,j$. For the dual $\Lambda$-basis 
${\tilde x}^{\prime*}_i\, (i=1,2,\dots, n)$ of 
$\tilde x^*_i\, (i=1,2,\dots, n)$ , let 
\begin{eqnarray*}
\tilde z_1^* &=& {\tilde x}^{\prime*}_1- b'_{11}(t)\tilde y^*_1,\\
\tilde z_2^* &=& {\tilde x}^{\prime*}_2- b_{12}(t)\tilde y^*_1-b'_{22}(t)\tilde y^*_2,
\\
&\dots& \\
\tilde z_n^* &=& {\tilde x}^{\prime*}_n- b_{1n}(t)\tilde y^*_1-b_{2n}(t)\tilde y^*_2-
\dots-b'_{nn}(t)\tilde y^*_n.
\end{eqnarray*}
Then we have 
\[\mbox{Int}^{\widetilde Y}_{\Lambda}(\tilde y^*_i,\tilde y^*_j)=0, \quad
\mbox{Int}^{\widetilde Y}_{\Lambda}(\tilde y^*_i,\tilde z^*_j)=\delta_{ij}, \quad
\mbox{Int}^{\widetilde Y}_{\Lambda}(\tilde z^*_i,\tilde z^*_j)
=\varepsilon^b_i\delta_{ij}\]
for all $i,j$. This $\Lambda$-basis creates 
a 3-sphere leaf $V^*$ for $Y^*$ and defines a winding degree 
$\lambda$ on $M$. 

The $\Lambda$-intersection numbers 
$\mbox{Int}^{\widetilde X}_{\Lambda}(\tilde y_i,\tilde y^*_j)$ and 
$\mbox{Int}^{\widetilde X}_{\Lambda}(\tilde y_i,\tilde z^*_j)$ 
are calculated as follows. 
\begin{eqnarray*}
\mbox{Int}^{\widetilde X}_{\Lambda}(\tilde y_i,\tilde y^*_j)&=& 
\mbox{Int}^{\widetilde X}_{\Lambda}({\tilde x'}_i+{\tilde x}^{\prime-}_i, 
\tilde x^-_j+\tilde x^*_j)= 
\mbox{Int}^{\widetilde X}_{\Lambda}({\tilde x}^{\prime-}_i, \tilde x^-_j)
=-\delta_{ij}.\\
\mbox{Int}^{\widetilde X}_{\Lambda}(\tilde y_i,\tilde z^*_j)
&=&
\mbox{Int}^{\widetilde X}_{\Lambda}(
{\tilde x'}_i+{\tilde x}^{\prime-}_i,
{\tilde x}^{\prime*}_j- \sum_{k'=1}^j 
\tilde b_{k'j}(t)(\tilde x^-_{k'}+\tilde x^*_{k'})\\
&=&
\mbox{Int}^{\widetilde X}_{\Lambda}({\tilde x}^{\prime-}_i,
-\sum_{k'=1}^j \tilde b_{k'j}(t) \tilde x^-_{k'})
= \tilde b_{ij}^0(t).
\end{eqnarray*}

The $\Lambda$-intersection numbers 
$\mbox{Int}^{\widetilde X}_{\Lambda}(\tilde z_i,\tilde y^*_j)$ and 
$\mbox{Int}^{\widetilde X}_{\Lambda}(\tilde z_i,\tilde z^*_j)$ 
are calculated as follows.
\begin{eqnarray*}
\mbox{Int}^{\widetilde X}_{\Lambda}(\tilde z_i,\tilde y^*_j)&=& 
\mbox{Int}^{\widetilde X}_{\Lambda}(\tilde x^-_i- 
\sum_{k=1}^{i}\tilde a_{ik}(t^{-1})({\tilde x'}_k
+{\tilde x}^{\prime-}_k), 
\tilde x^-_j+\tilde x^*_j)\\
&=& 
\mbox{Int}^{\widetilde X}_{\Lambda}(\tilde x^-_i- 
\sum_{k=1}^{i}\tilde a_{ik}(t^{-1}){\tilde x}^{\prime-}_k, 
\tilde x^-_j)
= -a_{ij}(t)+\tilde a^{0*}_{ij}(t).\\
\mbox{Int}^{\widetilde X}_{\Lambda}(\tilde z_i,\tilde z^*_j)
&=&
\mbox{Int}^{\widetilde X}_{\Lambda}(
\tilde x^-_i- \sum_{k=1}^i \tilde a_{ik}(t^{-1})({\tilde x'}_k
+{\tilde x}^{\prime-}_k),
{\tilde x}^{\prime*}_j- 
\sum_{k'=1}^j \tilde b_{k'j}(t)(\tilde x^-_{k'}+\tilde x^*_{k'})\\
&=&
\mbox{Int}^{\widetilde X}_{\Lambda}(
\tilde x^-_i- \sum_{k=1}^i 
\tilde a_{ik}(t^{-1}){\tilde x}^{\prime-}_k,
-\sum_{k'=1}^j \tilde b_{k'j}(t) \tilde x^-_{k'})\\
&=&
\tilde b_{ij}^0(t)
-\sum_{k=1}^{\min\{i,j\}}\tilde a_{ik}(t)\cdot\tilde b_{kj}(t)\\
&=&
\tilde b_{ij}^0(t)
-\sum_{k=1}^{\min\{i,j\}}\tilde a_{ik}(t)\cdot\tilde b_{kj}(t)
=\tilde b_{ij}^0(t)-c_{ij}(t).
\end{eqnarray*}
By examining these calculations (particularly, by noting that 
 $-a_{ij}(t)+\tilde a^{0*}_{ij}(t)$ means $-a_{ij}(t)\,(i<j)$, 
$-\varepsilon^a_i-a'_{ii}(t^{-1})$ or $0\, (i>j)$, and 
$\tilde b_{ij}^0(t)$ means $b_{ij}(t)\,(i<j)$, $b'_{ii}(t)$ 
or $0\, (i>j)$, it is seen that 
$\lambda\leqq \lambda(A)$. 
$\square$

\phantom{x}

\noindent{\bf 3. Proofs of Theorems~1.1 and 1.4} 

It is not always assumed that a closed 4-manifold is a 
smooth or piecewise-linear manifold, but smooth and piecewise-linear 
techniques can be used  for it because a punctured manifold of it is 
smoothable (see \cite{FQ}). 

The proof of Theorem~1.1 will be basically  analogous to the proof 
of \cite[Theorem~1.1]{KA2}. 
For the use of a minimal element, we provide the following lemma. 

\phantom{x}

\noindent{\bf Lemma~3.1.} For a TOP-split positive definite 
$Z^{\pi_1}$-manifold $M$, 
let $x\in H_2(M;Z)$ be a minimal element. For any element 
$x'\in H_2(M^{(2)};Z)$ with $p_*(x')=x$, the inequality 
$||x'||^2 \geqq ||x||^2$ holds. 

\phantom{x}

\noindent{\bf Proof.} Let $M=S^1\times S^3\#M_1$ for 
a simply connected closed 4-manifold $M_1$, and 
$M^{(2)}=S^1\times S^3\#M_{1,1}\# M_{1,2}$ for the two copies 
$M_{1,1}, M_{1,2}$ of $M_1$. The element $x$ is represented 
as a 2-cycle in $M_1$.  Let $x_1$  be an element of 
$H_2(M^{(2)};Z)$ represented as a 2-cycle of $M_{1,1}$ 
such that $p_*(x_1)=x$ and 
$||x_1||^2 =||x||^2$. 
Then there are elements $y_i\, (i=1,2)$ of $H_2(M^{(2)};Z)$ 
such that $y_i$ is represented as 2-cycles in $M_{1,i}$ and 
$x'=x_1+y_1+y_2$. The identity  
$-p_*(y_1)=p_*(y_2)\in H_2(M;Z)$ is obtained from 
$p_*(x')=x$, meaning that $y_2$ 
is nothing but 
the covering translation element of $-y_1$ in $H_2(M^{(2)};Z)$. 
Denoting $-p_*(y_1)=p_*(y_2)$ by $y$, we have 
\[||y_1||^2=||y_2||^2=||y||^2.\]
Since $p_*(x_1+y_1)=x-y$ and $||x_1+y_1||^2 =||x-y||^2$, 
we have $||x'||^2=||x-y||^2+||y||^2$. 
It is noted that this equality holds without the positive definiteness. 
Using that $M$ is positive definite and $x$ is a minimal element with 
$x=(x-y)+y$, we must have 
\[||x-y||^2\geqq ||x||^2 \quad{\mbox{or}}\quad 
||y||^2\geqq ||x||^2.\]
Thus, the inequality $||x'||^2\geqq ||x||^2$ holds. 
$\square$

\phantom{x}

The proof of Theorem~1.1 will be  done as follows:

\phantom{x}

\noindent{\bf Proof of Theorem~1.1.}

\noindent{\bf Proof of (0) \mbox{$\mathbf \to$} (2).}\,  
Let $M=S^1\times S^3\#M_1$ for 
a simply connected closed 4-manifold $M_1$. Then $\widetilde M$ is the connected 
sum of $R^1\times S^3$ and the infinite copies $M_{1j}\, (j=0,\pm1,\dots)$ of 
$M_1$. Then every minimal element $\tilde x_i$ is represented by a $2$-cycle 
in one copy  $M_{1j}$ after a $t$-power shift, so that the elements $\tilde x_i\, (i=1,2,\dots,n)$
are represented by the same copy $M_{10}$ after suitable $t$-power shifts of 
$\tilde x_i\, (i=1,2,\dots,n)$, showing (0) $\to$ (2). $\square$

\phantom{x}

\noindent{\bf Proof of (2) \mbox{$\mathbf \to$} (1).}\,  
This assertion is obvious since after suitable $t$-shifts 
of $\tilde x_i\, (i=1,2,\dots,n)$, the $\Lambda$-intersection numbers 
$\mbox{Int}^{\widetilde M}_{\Lambda}(\tilde x_i, \tilde x_j)$ belong to $Z$ 
and the elements $\tilde x^{(1)}_i\, (i=1,2,\dots,n)$ are $Z$-generators 
for $H_2(M;Z)$. $\square$

\phantom{x}

\noindent{\bf Proof of (1) \mbox{$\mathbf \to$} (0).}\,  
This assertion will mean that the conditions (0), (1) and (2) are 
mutually equivalent. 
Assume that the elements $\tilde x_i\, (i=1,2,\dots,n)$ in $H_2(\widetilde M;Z)$ 
with $a_{ij}=\mbox{Int}^{\widetilde M}_{\Lambda}(\tilde x_i, \tilde x_j)\in Z$  
for all $i,j$ 
induce $Z$-generators $\tilde x_i^{(1)}\, (i=1,2,\dots,n)$ for $H_2(M;Z)$.  
Then we have the $Z$-intersection number 
$\mbox{Int}^{M}(\tilde x_i^{(1)}, \tilde x_j^{(1)})=a_{ij}$ 
for all $i,j$. Let $H_2(M;Z)$ be a $Z$-free group of rank $s$.   
For a $Z$-basis $y_j\, (j=1,2,\dots,s)$ of $H_2(M;Z)$, every basis element $y_j$ is 
a $Z$-linear combination of the elements $\tilde x_i^{(1)}\, (i=1,2,\dots, n)$. 
Let $\tilde y_j $ be the element in $H_2(\widetilde M;Z)$ to be the 
$Z$-linear combination on $\tilde x_i\, (i=1,2,\dots, n)$ whose  $Z$-coefficients are the same 
 as the $Z$-coefficients in the $Z$-linear combination of $y_j$ on 
$\tilde x_i^{(1)}\, (i=1,2,\dots, n)$. 
Then the $\Lambda$-intersection matrix  
$\left(\mbox{Int}_{\Lambda}^{\widetilde M}(\tilde y_j,\tilde y_{j'})\right)$ 
is equal to 
the $Z$-intersection matrix $(\mbox{Int}^{M}(y_j,y_{j'}))$, 
which is a unimodular matrix of size $s$. Using that 
$H_2(\widetilde M;Z)$ is a $\Lambda$-free module of rank $s$, 
we see that the elements $\tilde y_j\, (j=1,2,\dots,s)$ form  
a $\Lambda$-basis for $H_2(\widetilde M;Z)$.  Then it is shown in  \cite[Corollary~3.4]{K1} 
that  $M$ is TOP-split, but reproved here for convenience.  
In fact,  for $Z^{\pi_1}$-manifold  $M_0=S^1\times S^3\# M_1$ with $M_1$  the simply connected 4-manifold 
obtained from $M$ by a $2$-handle surgery killing $\pi_1(M)\cong Z$,  
the $\Lambda$-intersection forms on $\widetilde M$ and $\widetilde M_0$  
are $\Lambda$-isomorphic. Since $M$ and $M_0$ have the same Kirby-Siebenmann obstruction, 
$M$ and $S^1\times S^3\# M_1$ are homeomorphic  by \cite{FQ}. 
Thus, $M$ is TOP-split, showing (1)$\to$ (0). 
$\square$

\phantom{x}

\noindent{\bf Proof of (0) \mbox{$\mathbf \to$} (3).}\,  
For $M=S^1\times S^3\#M_1$ 
and $\widetilde M=R^1\times S^3\#_{j=-\infty}^{+\infty} M_{1j}$ 
(already established), every minimal 
element $\widetilde x\in H_2(\widetilde M;Z)$ is represented by a $2$-cycle $c_x$
in the copy $M_{10}$ after a $t$-power shift, 
which is not null-homologous in $M_{10}$. Assume that the element $\widetilde x^{(m)} \in H_2(M^{(m)};Z)$ 
is  the sum $y+z$ for non-zero elements $y,z\in H_2(M^{(m)};Z)$. 
Let $y=[c'_y+c''_y], z=[c'_z+c''_z]\in H_2(M^{(m)};Z)$ where  
$c'_y, c'_z$ are 2-cycles in $M_{10}$ and $c''_y, c''_z$ are 2-cycles in the complement 
$\mbox{cl}(M^{(m)}\setminus M_{10})$. 
Since $\widetilde x=[c_x]=[c_y]+[c_z]\in H_2(\widetilde M;Z)$, we have 
$||\widetilde x||^2\leqq ||[c_y]||^2$ or 
$||\widetilde x||^2\leqq ||[c_z]||^2$ 
by the definition of minimal elements. 
This means that 
\[||\widetilde x^{(m)}||^2\leqq 
||[c_y]||^2\leqq ||y||^2 \quad \mbox{or} \quad ||\widetilde x^{(m)}||^2\leqq 
||[c_z]||^2\leqq ||z||^2.\]
Thus, the element $\widetilde x^{(m)}\in H_2(M^{(m)};Z)$ is a minimal element 
for every $m$.  Let $x'\in H_2(M^{(2m)};Z)$ be an element 
with $p_*(x')=\widetilde x^{(m)}$. Since $M^{(m)}$ is TOP-split, we see from 
Lemma~3.1 that 
$||x'||^2 \geqq ||\widetilde x^{(m)}||^2$, 
showing (0) $\to$ (3). $\square$

\phantom{x}

\noindent{\bf Proof of (3) \mbox{$\mathbf \to$} (4).}\,  
Let $\widetilde x_i\, (i=1,2,\dots,n)$  be minimal $\Lambda$-generators 
of $H_2(\widetilde M;Z)$. By (3), the elements 
\[t^j(\widetilde x_i^{(m)})\in H_2(M^{(m)};Z)\quad(i=1,2,\dots,n; j=0,1,\dots,m-1)\]
give desired minimal $Z$-generators of $H_2(M^{(m)};Z)$  without assumption 
that $m\geqq \lambda$, showing  (3) $\to$ (4).  $\square$

\phantom{x}

It will be directly shown in 
Proposition~3.2 that the inequality  
$||x'_i||^2 > ||x_i||^2-2$ in (4) 
is always equivalent to the inequality 
$||x'_i||^2\geqq ||x_i||^2$. 

\phantom{x}

\noindent{\bf Proof of (4) \mbox{$\mathbf \to$} (0).}\, 
By \cite{KA}, $M$ is TOP-split if and only if $M^{(m)}$ is TOP-split. 
The proof will be done by replacing $M^{(m)}$ with $m\geqq \lambda$ 
for a previously given  winding degree $\lambda$ on $M$ 
by $M$ with $1\geqq \lambda$.
Let $x_i\, (i=1,2,\dots,s)$ be minimal $Z$-generators of  
$H_2(M;Z)$ such that $||x'_i||^2 > ||x_i||^2-2$ for every element 
$x'_i\in H_2(M^{(2)};Z)$ with $p_*(x'_i)=x_i$ and every $i$. 
The winding degree $\lambda$ on $M$ is given by the number 
$d([uP]';V',X')$ on $X=\overline Q(u)\# Q(v)\# M=X'\circ X''$ as defined in \S~2.  
Let $P_k$ be the $k$th 2-sphere in the 2-sphere union $u P$ in $\overline Q(u)$.  
Let $\widetilde P_k$ be a connected lift of $P_k$ to $\widetilde X$. 
After a $t$-power shift, the homology class $[\widetilde P_k]'\in 
H_2(\widetilde X';Z)$ is written as 
\[[\widetilde P_k]'=[\widetilde PD_k]=c_{k,0} +c_{k,1} t,\]
where  $\widetilde PD_k$ is a connected lift of a closed connected oriented surface 
$PD_k=P^E_k\cup (-D_k)\subset X'$ constructed in \S~2 so that 
$P^E_k\subset \mbox{cl}(X'\setminus N')$ and 
$D_k\subset \partial N'\setminus \partial B'$, 
and $c_{k,i}\, (i=0,1)$ 
are homology classes in $H_2(\widetilde X';Z)$ represented by 
2-cycles in the fundamental region $X'_{V'}$. 

Let $\widetilde S^3$ be a fixed connected lift  of the 3-sphere leaf $S^3=V'$ of 
$X'$ to the infinite cyclic covering space $\widetilde X'$. 
Let $\widetilde L_{k,j}=t^{-j}(\widetilde PD_k\cap t^j \widetilde S^3)$ 
be an oriented link (possibly the empty link) in $\widetilde S^3$, and $L_{k,j}$ 
the projection of $\widetilde L_{k,j}$ into $S^3$. 
Represent the element $x_i$ by a closed connected oriented surface 
$F_i$ (embedded) in $M\setminus N$, which are disjoint from 
the surfaces $PD_k$. 
By modifying $F_i$ in a collar $S^3\times [-1,1]$ 
of the 3-sphere leaf $S^3$ in $X'$, assume that the intersection 
$K_i=F_i\cap S^3$ is an oriented trivial knot bounding a disk 
$\Delta_i$ in $S^3$ such that the system $\Delta_i\, (i=1,2,\dots, s)$ 
are mutually disjoint and meet $PD_k$ transversely in a finite number of points. 

For the preimage $\widetilde F_i$ of $F_i$ 
under the covering $\widetilde X'\to X'$, let 
$\widetilde K_i=\widetilde F_i\cap \widetilde S^3$, and $\widetilde \Delta_i$ the lift of 
$\Delta_i$ to $\widetilde S^3$ with $\partial \widetilde \Delta_i=\widetilde K_i$. 
We show that 
\[\ell_{i,k,j}={\mbox{Link}}^{\widetilde S^3}(\widetilde K_i,\widetilde L_{k,j})=0\] 
for all $i,j,k$. 
By Lemma~2.2, $\ell_{i,k,j}=0$ for 
$j\leqq 0$ or $j>1$ and all $i$ and $k$. 
Suppose $\ell_{i,k,1}$ is not zero for some $i$ and $k$. 

Construct an immersed disk $\Delta^*_i$ in $X$ with $\partial \Delta^*_i=K_i$ from
$\Delta_i$ by replacing a meridian disk of every point 
$L_{k,j}\cap\Delta_i$ for all $j$  in $\Delta_i$ with a disk disjointedly parallel 
to $P_k$ for every $k$.  
 
Let $G_i$ be a singular closed connected oriented surface in $X$ obtained 
from $F_i$ by cutting along $K_i$ and attaching two anti-oriented 
copies $\pm \Delta^*_i$ of $\Delta^*_i$. 
As an important note, by taking 
\[a_i=\sum_{j=-\infty}^{+\infty}\sum_{k=1}^u \ell_{i,,k,j}^2=\sum_{k=1}^u \ell_{i,k,1}^2>0,\]
the surface $\Delta^*_i$ has the $Z$-self-intersection number  $-a_i$ in $X$ with respect to 
the Seifert framing of $K_i$ in $S^3$. 
Let $x_{i,X}$ be the image of $x_i$ by the monomorphism $H_2(M;Z)\to H_2(X;Z)$. 
Then we have
\[x_{i,X}=[F_i]=[G_i]\in H_2(X;Z).\]
For a connected lift $G'_i$ of $G_i$ to the double covering space $X^{(2)}$ of $X$, 
let 
\[x'_{i,X}=[G'_i]\in H_2(X^{(2)};Z).\]
The square length of the element $x'_{i,X}$ is estimated as follows: 
\[(*)\qquad\qquad ||x'_{i,X}||^2=||x_{i,X}||^2-2a_i\leqq ||x_{i,X}||^2-2=||x_i||^2-2.\]
Since $\Delta^*_i \cap u P =\emptyset$, the surface $\Delta^*_i$ is regarded as a 
surface in the positive definite $Z^{\pi_1}$-manifold 
$M_Q=Q(v)\# M$ obtained from $X$ by blowing down on $\overline Q(u)$.
Let $x_{i,Q}$ be the  image of $x_i$ by the monomorphism $H_2(M;Z)\to H_2(M_Q;Z)$.  Then  we have 
\[x_{i,Q}=[F_i]=[G_i]\in H_2(M_Q;Z).\]
Since the surface $G'_i$ is in the connected summand $M_Q^{(2)}$ of $X^{(2)}$, let 
\[x'_{i,Q}=[G'_i]\in H_2(M_Q^{(2)};Z),\]
which is sent to $x_{i,Q}$ under the double covering projection homomorphism  $H_2(M_Q^{(2)};Z)\to H_2(M_Q;Z)$. 
Then the  inequality $(*)$  is equivalent to  the inequality 
\[(**)\qquad\qquad ||[x'_{i,Q}||^2\leqq ||x_{i,Q}||^2-2=||x_i||^2-2.\]
By using that  $H_2(M^{(2)};Z)$ is an orthogonal summand of $H_2(M^{(2)}_Q;Z)$,
let $x'_i$ be the image of  $x'_{i,Q}$ by the orthogonal summand  projection 
$H_2(M^{(2)}_Q;Z)\to H_2(M^{(2)};Z)$.  
Since $M^{(2)}_Q$ is positive definite, we have 
\[||[x'||^2\leqq ||x'_{i,Q}||^2,\quad
\mbox{so that}\quad  ||x'_i||^2\leqq ||x_i||^2-2.\] 
This  contradicts the inequality 
$||x'_i||^2> ||x_i||^2-2$ 
given by the assumption of (4),  
because the double covering projection homomorphism  $p_*: H_2(M^{(2)};Z)\to H_2(M;Z)$ 
sends $x'_i$   to $x_i$.  Thus, we have
\[\ell_{i,k,j}=\mbox{Link}^{\widetilde S^3}(\widetilde K_i,\widetilde L_{j,k})=0\] 
in the 3-sphere $\widetilde S^3$ for all $i,j,k$. 
This   means that 
all the $Z$-intersection numbers (containing the 
$Z$-self-intersection numbers) on any connected surfaces 
lifting the surfaces $\Delta^*_i (i=1,2,\dots,n)$   to 
the infinite cyclic covering space $\widetilde X$ have $0$ 
with respect to the lifted framings of 
the Seifert framings of $K_i  (i=1,2,\dots,n)$ in $S^3$. 
Let $\widetilde G_i$  be a  connected lift of the singular 
closed oriented surface $G_i$ in $X$ (already constructed above) to $\widetilde X$. 
Let 
\[\tilde x_{i,X}=[\widetilde G_i]\in H_2(\widetilde X;Z)\, (i=1,2,\dots, n). \]
Then,  from construction,  we see that the $\Lambda$-intersection numbers  
$\mbox{Int}_{\Lambda}^{\widetilde X}(\tilde x_{i,X},\tilde x_{i',X})$ 
are integers for all $i,i'$.
Since the surfaces $\widetilde G_i\, (i=1,2,\dots, n)$ belong to the infinite 
cyclic covering space $\widetilde M_Q$ of $M_Q$, 
let 
\[\tilde x_{i,Q}=[\widetilde G_i]\in H_2(\widetilde M_Q;Z)\, (i=1,2,\dots, n). \]
Note that the element  $\tilde x_{i,Q}$ is sent to $x_{i,Q}$ under the covering 
projection homomorphism $H_2(\widetilde M_Q;Z)\to H_2(M_Q;Z)$ and 
\[\mbox{Int}_{\Lambda}^{\widetilde X}(\tilde x_{i,X},\tilde x_{i',X})
=\mbox{Int}_{\Lambda}^{\widetilde M_Q}(\tilde x_{i,,Q},\tilde x_{i',Q})=
\mbox{Int}^{M_Q}(x_{i,Q}, x_{i',Q})  =\mbox{Int}^{M}(x_i, x_{i'}) \in Z\]
for all $i,i'$. 
Then, by an argument of  (1) $\to$ (0), there is  a free $\Lambda$-submodule  $\widetilde H$ 
of $H_2(\widetilde M_Q;Z)$ with a $\Lambda$-basis  $\tilde y_i\,(i=1,2,\dots, s^*)$ such that 

\phantom{x}

\noindent(i) \, every basis element $\tilde y_i$ is a $Z$-linear combination of  $\tilde x_{i,,Q}\, (i=1,2,\dots,s)$ 

\noindent(ii) \, every element $\tilde x_{i,,Q}$ is written as a $Z$-linear combination of $\tilde y_i\,(i=1,2,\dots, s^*)$.

\noindent(iii) \, the $\Lambda$-intersection matrix with respect to $\tilde y_i\,(i=1,2,\dots, s^*)$ is an integral matrix 
with determinant $+1$. 

\phantom{x}

We claim that there are minimal elements $\tilde z_i\, (i=1,2,\dots,s)$ in $\widetilde H$ 
such that  $\mbox{Int}_{\Lambda}^{\widetilde M_Q}(\tilde z_i,\tilde z_{i'})\in Z$ for all $i, i'$ and,  
by the covering projection 
homomorphism $H_2(\widetilde M_Q;Z)\to H_2(M_Q;Z)$, 
the elements $\tilde z_i\, (i=1,2,\dots,s)$
are sent to the minimal elements $x_{i,Q}\,  (i=1,2,\dots,s)$, respectively.  
In fact, by an argument on an orthogonal complement of the $\Lambda$-intersection form, there is an orthogonal  splitting 
\[H_2(\widetilde M_Q;Z)=\widetilde H\oplus \widetilde H'\]
with respect to the $\Lambda$-intersection form $\mbox{Int}_{\Lambda}$ on $H_2(\widetilde M_Q;Z)$. 
By Freedman-Quinn construction in  \cite{FQ},  there is a circle union splitting  
$Y\circ Y'$ of $M_Q$  for some  positive definite $Z^{\pi_1}$-manifolds  $Y, Y'$ such that 
the $\Lambda$-intersection form on $H_2(\widetilde Y;Z)$ is isomorphic to the restriction of the $\Lambda$-intersection form 
on $H_2(\widetilde M_Q;Z)$ to $\widetilde H$. By (1) $\to$ (0),  $Y$ is a TOP-split.  Let $Y_V$ be a fundamental region 
of $Y$ splitting along a 3-sphere leaf $V$. Let  $\tilde z_i\, (i=1,2,\dots,s)$  be  $\Lambda$-generators  
of $H_2(\widetilde Y;Z)=\widetilde H$ represented by 2-cycles in $Y_V$ and sent respectively to the minimal elements
$x_{i,Q}\,  (i=1,2,\dots,s)$.  Because $H_2(Y_V;Z)$  is isomorphic to the orthogonal summand $H_2(M;Z)$ of 
$H_2(\widetilde M_Q;Z)$ by the covering projection homomorphism 
$H_2(\widetilde M_Q;Z)\to H_2(M_Q;Z)$ and $H_2(Y_V;Z)$ is an orthogonal summand of 
$H_2(\widetilde Y;Z)=\widetilde H$, it follows that the elements $\tilde z_i\, (i=1,2,\dots,s)$ are minimal elements in 
$\widetilde H$, as desired. 
 
Using  that $\widetilde H$ is an orthogonal summand of $H_2(\widetilde M_Q;Z)$, we see that the elements 
$\tilde z_i\, (i=1,2,\dots,s)$ are minimal elements in $H_2(\widetilde M_Q;Z)$. 
Every minimal element of $H_2(CP^2;Z)$ must be a generator, so that  the minimal element 
$\tilde z_i\in H_2(\widetilde M_Q;Z)$ is represented by a 2-cycle in $\widetilde M$ because $\tilde z_i$ is 
sent to $x_{i,Q}$. This implies  that there are elements
$\tilde x_i\, (i=1,2,\dots,s)$ in  $H_2(\widetilde M;Z)$ coming from $\tilde z_i\, (i=1,2,\dots,s)$ such that 
the elements $\tilde x_i^{(1)}, (i=1,2,\dots,s)$ are equal to the $Z$-generators $x_i\, (i=1,2,\dots,s)$ 
of  $H_2(M;Z)$ coming from  $x_{i,Q}\,  (i=1,2,\dots,s)$ and the $\Lambda$-intersection numbers 
$\mbox{Int}^{\widetilde M}_{\Lambda}(\tilde x_i, \tilde x_{i'})$ 
are integers for all $i,i'$. 
By (1) $\to$ (0), $M$ is TOP-split, showing (4) $\to$ (0). $\square$

\phantom{x}

\noindent{\bf Proof of (0) \mbox{$\mathbf \to$} (5).}\,  If $M$ is TOP-split, then there is a 
winding degree $\lambda$ on $M$ with  $\lambda=0$ by definition. $\square$

\phantom{x}

\noindent{\bf Proof of (5) \mbox{$\mathbf \to$} (0).}\, 
Assume that $\lambda_{\mbox{min}}=\lambda=d([uP]';V',X')=0$ on 
$X=\overline Q(u)\# Q(v)\# M=X'\circ X''$ as defined in \S~2. 
The proof is almost similar to  the proof of the assertion (4) $\to$ (0). 
Let $P_k$ be the $k$th 2-sphere in the 2-sphere union $u P$ in $\overline Q(u)$.  
Let $\widetilde P_k$ be a connected lift of $P_k$ to $\widetilde X$. 
After a $t$-power shift, the homology class $[\widetilde P_k]'\in 
H_2(\widetilde X';Z)$ is written as 
\[[\widetilde P_k]'=[\widetilde PD_k]=c_{k,0},\]
where  $\widetilde PD_k$ is a connected lift of a closed connected oriented surface 
$PD_k=P^E_k\cup (-D_k)\subset X'$ constructed in \S~2 so that 
$P^E_k\subset \mbox{cl}(X'\setminus N')$ and 
$D_k\subset \partial N'\setminus \partial B'$, 
and $c_{k,0}$ is a homology class in $H_2(\widetilde X';Z)$ represented by 
a 2-cycle in the fundamental region $X'_{V'}$. 

Let $\widetilde S^3$ be a fixed connected lift  of the 3-sphere leaf $S^3=V'$ of 
$X'$ to the infinite cyclic covering space $\widetilde X'$. 
Let $\widetilde L_{k,j}=t^{-j}(\widetilde PD_k\cap t^j \widetilde S^3)$ 
be an oriented link (possibly the empty link) in $\widetilde S^3$, and $L_{k,j}$ 
the projection of $\widetilde L_{k,j}$ into $S^3$.

Let  $x_i\,(i=1,2,\dots,s)$ be any minimal $Z$-generators of  $H_2(M;Z)$.
Represent the element $x_i$ by a closed connected oriented surface 
$F_i$ (embedded) in $M\setminus N$, which are disjoint from 
the surfaces $PD_k$. 
By modifying $F_i$ in a collar $S^3\times [-1,1]$ 
of the 3-sphere leaf $S^3$ in $X'$, assume that the intersection 
$K_i=F_i\cap S^3$ is an oriented trivial knot bounding a disk 
$\Delta_i$ in $S^3$ such that the system $\Delta_i\, (i=1,2,\dots, s)$ 
are mutually disjoint and meet $PD_k$ transversely in a finite number of points. 

For the preimage $\widetilde F_i$ of $F_i$ 
under the covering $\widetilde X'\to X'$, let 
$\widetilde K_i=\widetilde F_i\cap \widetilde S^3$, and $\widetilde \Delta_i$ the lift of 
$\Delta_i$ to $\widetilde S^3$ with $\partial \widetilde \Delta_i=\widetilde K_i$. 
Then, by Lemma~2.2,  we have 
\[\ell_{i,k,j}={\mbox{Link}}^{\widetilde S^3}(\widetilde K_i,\widetilde L_{k,j})=0\] 
for all $i,j,k$.  

The rest of the proof is completely the same as the proof  of the assertion (4) $\to$ (0). 
This shows the assertion (5) $\to$ (0).  $\square$

\phantom{x}

This completes the proof of Theorem~1.1. $\square$

\phantom{x}

As an additional note to Theorem~1.1, the following proposition clarifies 
a reason why 
the square length $||x||^2$ in (2) may be replaced by $||x_i||^2-2$ in (3). 

\phantom{x}

\noindent{\bf Proposition~3.2.} For every $Z^{\pi_1}$-manifold $M$, 
any elements $x\in H_2(M;Z)$ and $x'\in H_2(M^{(2)};Z)$ 
with $p_*(x')=x$ have the congruence 
$||x'||^2 \equiv ||x||^2 \pmod{2}$. 

\phantom{x}

\noindent{\bf Proof.} For a {\it smooth} $Z^{\pi_1}$-manifold $M$, 
this property follows from the fact that 
every covering preserves the second Wu class. Let $M_0$ be a TOP-split 
$Z^{\pi_1}$-manifold. Then this property on $M_0$ holds as it is seen 
from a calculation of the $Z$-self-intersection numbers in the proof of Lemma~3.1. 
For a general $Z^{\pi_1}$-manifold $M$, there is  
a $Z$-homology cobordism $W$  from $M$ to a TOP-split  $Z^{\pi_1}$-manifold $M_0$,  
which is seen from the proof of \cite[Theorem~1.1]{K0} although this theorem itself contains 
a serious error (cf. \cite{KK, K1, K2}). Then the double covering space 
$W^{(2)}$ gives a $Z/2Z$-homology cobordism  from $M^{(2)}$ to 
$M_0^{(2)}$. The pair $(W^{(2)}, W)$ sends the pair of elements $x',x$ with 
$p_*(x')=x$ to a pair of elements 
$x'_0\in H_2(M^{(2)}_0;Z), x_0\in H_2(M^{(2)}_0;Z)$ 
with $p_*(x'_0)=x_0$ up to $2$ times elements. Thus, the congruence 
$||x'||^2\equiv ||x||^2 \pmod{2}$ 
is obtained. $\square$

\phantom{x}

The following corollary is slightly stronger than Corollary~1.2 and obtained from  Theorem~1.1 (4) 
since, as noted in \S~1, if there are  
$Z$-generators $x_i\, (i=1,2,\dots, n)$ of $H_2(M^{(m)};Z)$ with 
$||x_i||^2\leqq 2$ for all $i$, then there are minimal $Z$-generators $y_j\, (j=1,2,\dots, s)$ of  
$H_2(M^{(m)};Z)$ with $||y_j||^2\leqq 2$ for all $j$ and,  for every element $y'_j\in H_2(M^{(2m)};Z)$ 
with $p_*(y'_j)=y_j$, the inequality $||y'_j||^2 > ||y_j||^2-2$ holds.  

\phantom{x}

\noindent{\bf Corollary~3.3.}  A positive definite $Z^{\pi_1}$-manifold $M$ is 
TOP-split if  for any previously given winding degree $\lambda$ on $M$, 
there is an $m\geqq \lambda$ for which there are $Z$-generators 
$x_i\,(i=1,2,\dots, n)$ of $H_2(M^{(m)};Z)$ such that  $||x_i||^2\leqq 2$  
for all $i$. 

\phantom{x}

The proof of Theorem~1.4 will be  done as follows:

\phantom{x}

\noindent{\bf Proof of Theorem~1.4.} 
Assume that $M$ is TOP-split. Then there 
is a $\Lambda$-basis  $\tilde e_i\in H_2(\widetilde M;Z)\, (i=1,2,\dots,n)$ with 
$\mbox{Int}^{\widetilde M}_{\Lambda}(\tilde e_i,\tilde e_j)=\delta_{ij}$ 
for all $i,j$. For a $\Lambda$-basis $\tilde x_i\in H_2(\widetilde M;Z)\, 
(i=1,2,\dots,n)$, let $\tilde x_i=\sum_{k=1}^s a_{ik}(t)\tilde e_k$. Then 
we have 
\[\mbox{Int}^{\widetilde M}_{\Lambda}(\tilde x_i,\tilde x_j)
=\sum_{k=1}^n a_{ik}(t^{-1})a_{jk}(t),\]
showing (0) $\to$ (1). Assume (1). Let $P$ be the matrix of size $n$ whose $(i,j)$ 
entry is $a_{ij}(t)$. Then the $\Lambda$-intersection matrix $S$ 
whose $(i,j)$ entry is 
$\mbox{Int}^{\widetilde M}_{\Lambda}(\tilde x_i,\tilde x_j)$ 
is given by $P \overline P^{T}$. Since the determinant $\mbox{det}S=1$, the matrix 
$P$ is non-singular in $\Lambda$, so that $P^{-1}S (\overline P^{T})^{-1}$ is 
the identity matrix and  
$M$ is TOP-split  by for example Corollary~1.3, showing (1) $\to$ (0). 

The assertion (0) $\to$ (2) is direct. The converse (2) $\to$ (0) 
is obtained from Corollary~1.3 because there are  finitely many 
minimal $\Lambda$-generators of $H_2(\widetilde M;Z)$. 

(1) implies (3${}_1$) since every element of $H_2(\widetilde M;Z)$ is 
a $\Lambda$-linear combination of any given $\Lambda$-basis of 
$H_2(\widetilde M;Z)$. 
Also, (2) implies (3${}_2$) since there are
$\Lambda$-linearly independent minimal elements 
$\tilde y_i\, (i=1, 2,\dots, n)$. Thus, we have (0) $\to$ (3). 
To show (3) $\to$ (0), assume (3). 
Let $\tilde x_i\, (i=1,2,\dots,n-1)$  be mutually 
distinct elements of $H_2(\widetilde M;Z)$ 
up to multiplications of the units of $\Lambda$ such that 
the square length $||\tilde x_i||^2=1$ for all $i$. 
It is noted that 
if an element $f(t)\in\Lambda$ is written as $\sum_{i=p}^q c_i t^i\, 
(c_i\in Z)$, then the constant term $c$ of the product 
$f(t)f(t^{-1})$ is written as the sum $\sum_{i=p}^q c_i^2$, so that 
according to whether $c=0$ or $1$, we have $f(t)f(t^{-1})=0$ or $1$, 
respectively. This means that the $\Lambda$-square length 
$||\tilde x_i||^2_{\Lambda}=1$  
for all $i$. For the $\Lambda$-submodule $\Lambda[\tilde x_1]$ of 
$H_2(\widetilde M;Z)$ generated by $\tilde x_1$, 
let $H_2(\widetilde M)=\Lambda[\tilde x_1]\oplus H'$
be the orthogonal splitting
with respect to the $\Lambda$-intersection form 
$\mbox{Int}^{\widetilde M}_{\Lambda}$ on $H_2(\widetilde M;Z)$. 
Then $\tilde x_2$ is written as $a_1(t) \tilde x_1+\tilde x'_2$ 
for an element $a_1(t)\in\Lambda$  and an element $\tilde x'_2$ in $H'$. 
The square length 
$||\tilde x_2||^2$ is computed from the constant 
terms of the following identities:
\[||\tilde x_2||^2_{\Lambda}=a_1(t)a_1(t^{-1})+
||\tilde x'_2||^2_{\Lambda}= 1.\]
If $\tilde x'_2\ne 0$, then 
the square length  
$||\tilde x'_2||^2>0$, 
so that $a_1(t)=0$ and $\tilde x_2=\tilde x'_2$. If $\tilde x'_2=0$, then
$a_1(t)=\pm t^k$ for some integer $k$, so that $\tilde x_1$ and $\tilde x_2$ 
are equal up to the 
multiplication of a unit of $\Lambda$, which is  a contradiction. Thus, 
the identity 
$\tilde x_2=\tilde x_2'\in H'$ is obtained. 
By a similar argument, the elements $\tilde x_i\, (i=2,3,\dots,n-1)$ are seen 
to be in $H'$. 
By an inductive argument, it is shown that the elements 
$\tilde x_i\, (i=1,2,\dots,n-1)$ have the identities 
$\mbox{Int}^{\widetilde M}_{\Lambda}(\tilde x_i,\tilde x_j)=\delta_{ij}$ 
for all $i,j$ and hence form a $\Lambda$-basis of  the  $\Lambda$-submodule 
$\Lambda[\tilde x_1,\tilde x_2,\dots, \tilde x_{n-1}]$ of $H_2(\widetilde M;Z)$ generated 
by $\tilde x_i\, (i=1,2,\dots,n-1)$.  Let $H'_1$ be the orthogonal complement 
of $\Lambda[\tilde x_1,\tilde x_2,\dots, \tilde x_{n-1}]$  in 
$H_2(\widetilde M;Z)$, which is $\Lambda$-isomorphic to $\Lambda$ with $\Lambda$-intersection 
matrix $(1)$. Thus, there is a $\Lambda$-basis $\tilde x_i\, (i=1,2,\dots,n)$ 
of $H_2(\widetilde M)$ with  identities 
$\mbox{Int}^{\widetilde M}_{\Lambda}(\tilde x_i,\tilde x_j)=\delta_{ij}$ 
for all $i,j$. 
For example by Corollary~1.3, $M$ is TOP-split, showing (3) $\to$ (0). 
$\square$

\phantom{x}

\noindent{\bf 4. Proofs  of  Theorem~1.5, Corollary~1.6 and  Observations~1.7 and 1.8}

The proof of Theorem~1.5 will be  done as follows: 

\phantom{x}

\noindent{\bf Proof of Theorem~1.5.} 
The square length $||\tilde x_i||^2$ is $3$ for 
$i=1,2$ and $2$ for $i=3,4$.   
If $\tilde x_1$ is not minimal, then 
there are $\Lambda$-generators 
$\tilde x'_{i'}\, (i'=1,2,\dots,k)$ of $H_2(\widetilde M_L;Z)$ with 
$||\tilde x'_{i'}||^2 \leqq 2$ for all $i'$ 
because 
$||\tilde x_1-\tilde x_2||^2=2$. 
Then by  Corollary~1.3, 
the $Z^{\pi_1}$-manifold $M_L$ must be TOP-split, which contradicts 
that $M_L$ is not TOP-split. 
Hence $\tilde x_1$ is a  minimal element. By a similar argument, 
$\tilde x_2$ is also a  minimal element. 
If $\tilde x_3$ is not minimal, then 
$\tilde x_3$ is the sum of two elements $\tilde x'_3, \tilde x''_3$ with  
square lengths 
\[||\tilde x'_3||^2=||\tilde x''_3||^2=1.\]
The elements $t^i \tilde x_1\, (i\in Z)$ must belong to the 
$Z$-orthogonal complement of the infinite cyclic group generated by  $\tilde x'_3$ 
since they are minimal elements with 
square length $3$. Hence 
$\mbox{Int}^{\widetilde M_L}_{\Lambda}(\tilde x_1, \tilde x'_3)=0$. Similarly, 
$\mbox{Int}^{\widetilde M_L}_{\Lambda}(\tilde x_1, \tilde x''_3)=0$, so that 
$\mbox{Int}^{\widetilde M_L}_{\Lambda}(\tilde x_1, \tilde x_3)=0$ 
contradicting that 
$\mbox{Int}^{\widetilde M_L}_{\Lambda}(\tilde x_1, \tilde x_3)\ne 0$. Hence 
$\tilde x_3$ is a minimal element. By a similar method, $\tilde x_4$ is 
also a minimal element, showing (1). To see (2), suppose $\mbox{det}V_G=1$. 
The set $B=\{ t^i \tilde x_k |\, i\in Z, k=1,2,3,4\}$ forms a $Z$-basis for 
$H_2(\widetilde M_L;Z)$ consisting of minimal elements.  
By a property of a minimal element, every element of $B$ belongs to either $G$ or 
the $Z$-orthogonal complement 
$G^{\bot}$ of $G$ constructed by using $\mbox{det}V_G=1$. However, 
it is impossible 
because any two elements 
of $B$ are connected by a sequence $v_j\in B\, (j=1,2,\dots, s)$ 
with $\mbox{Int}^{\widetilde M_L}(v_j, v_{j+1})\ne 0$ for every $j$.  
Thus, $\mbox{det}V_G>1$, showing (2). 
$\square$

\phantom{x}

The proof of Corollary~1.6 will be  done as follows: 

\phantom{x}

\noindent{\bf Proof of Corollary~1.6.} Since 
$m\geqq 3\geqq e({\tilde x_i})+1\,\, (i=1,2,3,4)$, 
we have the square length 
\[||\tilde x_i^{(m)}||^2=||\tilde x_i||^2\quad (i=1,2,3,4)\]
which is $3$ for $i=1,2$ and $2$ for $i=3,4$.  
A winding degree $\lambda$ on $M_L$ with $\lambda\leqq 5$ 
is found from Lemma~2.3  by considering 
the following matrices:
\[
A=\left(
\begin{array}{cccc}
1+f+f^2 & 1  & 1+f& f\\
1       & 2  & 1  & -1\\
1+f     & 1  & 2  & 0\\
f       & -1 & 0  & 2
\end{array}
\right),
\]

\[
A^{-1}=\left(
\begin{array}{cccc}
4     & -2  & -1-2f  & -1-2f\\
-2    & 2   & f      & 1+f\\
-1-2f & f   & 1+f+f^2& f+f^2\\
-1-2f & 1+f & f+f^2  & 1+f+f^2
\end{array} 
\right),
\]
where the matrix $A$ is obtained as the $\Lambda$-intersection matrix of 
$H_2(\widetilde M_L;Z)$ given by the $\Lambda$-basis 
$\tilde x_1, \tilde x_1-\tilde x_2, \tilde x_3, \tilde x_4$, 
and $A^{-1}$ is obtained as the inverse matrix of $A$. 

The numbers $\max\lambda(A)=3$ and $\min\lambda(A)=-2$ are determined by actual 
calculations, where $\max\lambda(A)=3$ is attained only by the Laurent polynomials 
\[c_{13}(t)=c_{14}(t)=(1+t+t^2)(1+2f),\] 
and $\min\lambda(A)=-2$ is attained  
only by the following Laurent polynomials:
\begin{eqnarray*}
b_{34}(t)&=& f+f^2,\quad a'_{11}(t^{-1})=1+t^{-1}+t^{-2}, \\
c_{33}(t)&=&-(1+2f+2f^2)+(1+t+t^2),\quad c_{34}(t)=-f-f^2, \\
\quad c_{43}(t)&=&-2f-2f^2,\quad c_{44}(t)=-1-2f-2f^2.
\end{eqnarray*}
Thus, 
\[\lambda(A)=\max\lambda(A)-\min\lambda(A)=3-(-2)=5.\]

For $m\geqq 5\geqq \lambda$, it is seen that 
$\tilde x_i^{(m)}$ is minimal for $i=1,2$ by 
a similar consideration of the proof of Theorem~1.5 (1) using 
Corollary~3.3 instead of Corollary~1.3 since
\[||\tilde x_1^{(m)}-\tilde x_2^{(m)}||^2=||\tilde x_1-\tilde x_2||^2 =2.\] 
Also, by 
a similar consideration of the proof of Theorem~1.5 (1), 
$\tilde x_i^{(m)}$ is minimal for $i=3,4$, showing (1) for $m\geqq 5$. 
To see (2) for $m\geqq 5$, let $G$ be a proper $Z$-free subgroup of 
$H_2(M_L^{(m)};Z)$. Since the $Z$-basis 
$t^j \tilde x_i^{(m)}\,(i=1,2,3,4;j=0,1,2,\dots,m-1)$ of 
$H_2(M_L^{(m)};Z)$ are minimal, 
a similar consideration of the proof of Theorem~1.5 (2) shows 
that the $Z$-intersection form on $G$ has the determinant greater 
than $1$, showing (2) for $m\geqq 5$. The assertion 
(2) for $m=3, 4$ was given in \cite{FHMT}. By assuming it, the assertion (1) for $m=3,4$ is shown as follows: 

The elements $\tilde x_i^{(m)}$ for $i=3,4$ must be minimal 
because $||\tilde x_i^{(m)}||^2=2$ and every element of 
the square length $1$ generates an orthogonal summand $Z$. 
Suppose $\tilde x_1^{(m)}$ is written as a sum $x'_1+x''_1$ with 
$||x'_1||^2<3$ and $||x''_1||^2<3$. Then 
$||x'_1||^2=||x''_1||^2=2$ which contradicts 
 $||\tilde x_1^{(m)}||=3$. Hence $\tilde x_1^{(m)}$ is minimal. 
Similarly,  $\tilde x_2^{(m)}$ is shown to be minimal. 
$\square$

\phantom{x}

The proof of Observation~1.7 will be  done as follows: 

\phantom{x}

\noindent{\bf Proof of Observation~1.7.} 
Suppose the $3$-manifold $nB$ bounds a smooth compact oriented $4$-manifold 
$W$ with $H_2(W;Q)=0$. 
Since $H_2(W,nB;Q)=H^2(W;Q)=0$ by Poincar{\'e} duality, the natural 
map $H_1(nB;Z)\to H_1(W;Z)/(\mbox{torsion})$ is injective. 
Then there is an epimorphism $\gamma:H_1(W;Z)\to Z$ such that 
the infinite cyclic covering $\widetilde W\to W$ belonging to $\gamma$ 
lifts at least one component of $nB$ non-trivially. Let $\Lambda_Q=Q[t,t^{-1}]$. 
Since $H_2(W;Q)=0$, the $\Lambda_Q$-module $H_2(\widetilde W;Q)$ is a finitely generated 
torsion $\Lambda_Q$-module for  whose rational 
Alexander polynomial $A(t)\in \Lambda_Q$ has $A(1)\ne 0$. Since 
there are infinitely many positive integer $m$ such that $t^m-1$ is coprime with 
$A(t)$ in $\Lambda_Q$, 
the $m$-fold cyclic covering space $W^{(m)}$ of $W$ belonging to 
the epimorphism $H_1(W;Z)\stackrel{\gamma}{\to} Z\to Z/mZ$  has 
\[H_2(W^{(m)};Q)=H_2(\widetilde W;Q)/(t^m-1)H_2(\widetilde W;Q)=0\] 
for any such $m$.  
The boundary $\partial W^{(m)}$ of the 4-manifold $W^{(m)}$ is a disjoint union of  
$m_i$-fold cyclic covering spaces $B^{(m_i)}\, (i=1,2,\dots,s)$ of $B$ which are 
still homology handles. It is noted that 
some $m_i$ may be taken sufficiently large when $m$ is taken sufficiently large. 
Then a smooth compact oriented  $4$-manifold $Y$ with $H_1(Y;Q)=H_2(Y;Q)=0$ is 
constructed from $W^{(m)}$ by $2$-handle 
surgeries on  $W^{(m)}$ killing $H_1(W^{(m)};Q)$. 
Let $Y^*$ be the smooth closed oriented 4-manifold obtained from 
$-Y$ by attaching $E^{(m_i)}\, (i=1,2,\dots,s)$ along $B^{(m_i)}\, 
(i=1,2,\dots,s)$, where $E^{(m_i)}$ denotes the $m_i$-fold cyclic covering space 
of a smooth compact connected oriented  4-manifold $E$ with boundary $B$ explained 
in the introduction. 
By Corollary~1.6, for a large $m$ the smooth $4$-manifold $Y^*$ with $H_1(Y^*;Q)=0$ 
must have  a positive definite 
non-standard $Z$-intersection form on $H_2(Y;Z)/(\mbox{torsion})$, which contradicts 
\cite{D2}. Thus, the $3$-manifold $nB$ cannot bound any smooth compact oriented 
$4$-manifold $W$ with $H_2(W;Q)=0$. 
$\square$

\phantom{x}

The proof of Observation~1.8 will be  done as follows: 

\phantom{x}

\noindent{\bf Proof of Observation~1.8.}
Suppose that 
the knot $K$ in $S^3$ is made a trivial knot by changing 
$r$ positive crossings 
into negative crossings. 
Then the $(3,1)$-manifold pair $(-S^3, K)$, where $-S^3$ denotes $S^3$ with 
orientation reversed, bounds a $(4,2)$-manifold pair $(Y, D)$, 
where $Y=D^4\# Q(r)$ for the 4-disk $D^4$ and  $D$ is a smooth disk 
with $[D]=0$ in $H_2(Y,-S^3;Z)$. In this construction, the inclusion 
$-S^3\setminus k \subset Y\setminus D$ is assumed to induce 
an isomorphism $H_1(-S^3\setminus k;Z)\to H_1(Y\setminus D;Z)$ 
and an epimorphism 
$\pi_1(-S^3\setminus k)\to \pi_1(Y\setminus D)$. 
Let $W=\mbox{cl}(Y\setminus N)$ for a tubular neighborhood $N$ of $D$ in $Y$.  
Then $W$ is a smooth positive definite  compact 4-manifold with boundary $-B$ 
such that the inclusion $-B\subset W$ induces an isomorphism 
$H_1(-B;Z)\to H_1(W;Z)$ and an epimorphism $\pi_1(-B)\to \pi_1(W)$. 
Let $X$ be a smooth positive definite $Z^{\pi_1}$-manifold obtained 
$E$ by attaching $W$ along $B$. By Corollary~1.6, the 
$Z$-intersection form on $H_2(X^{(m)};Z)$ for the $m$-fold cyclic covering space 
$X^{(m)}$ of $X$ for $m\geqq 3$ which is a smooth $Z^{\pi_1}$-manifold  must have  
a positive definite non-standard form as an orthogonal summand, which contradicts 
\cite{D2}. Thus, the knot $k$ cannot be made a trivial knot by changing 
positive crossings into negative crossings. 
$\square$

\phantom{x}

\noindent{\bf Acknowledgement.} This work was supported by JSPS KAKENHI Grant 
Number 26287013.  The author thanks to a referee  for a very careful reading.

\end{document}